%% file: main.tex
\documentclass[12pt]{amsart}

\usepackage{amsmath,amssymb}
\usepackage{graphicx}
\usepackage{enumitem}
\usepackage{hyperref}
\usepackage{xcolor}
\usepackage[backend=biber, style=numeric, url=true]{biblatex}
\usepackage{url}
\usepackage{imakeidx}
\makeindex[title=Some Appearances of Frame Terminology and Related Concepts, intoc, columns=1]

\newcommand{\comment}[1]{}

\newcommand\ip[1]{\left\langle#1\right\rangle}
\newcommand\set[1]{\left\{#1\right\}}
\newcommand\norm[1]{\left\lVert#1\right\rVert}
\newcommand\p[1]{\left(#1\right)}

\newcommand{\abs}[1]{\left\lvert#1\right\rvert}

\newcommand{\floor}[1]{\left\lfloor#1\right\rfloor}

\newcommand{\keepmath}[1]{\hfil\penalty200\hfilneg$\relpenalty=10000 \binoppenalty=10000 #1$}

\newtheorem{thm}{Theorem}[section]
\newtheorem*{thm*}{Theorem}
\newtheorem{prop}[thm]{Proposition}
\newtheorem{lem}[thm]{Lemma}
\newtheorem*{lem*}{Lemma}
\newtheorem{cor}[thm]{Corollary}
\newtheorem*{cor*}{Corollary}
\theoremstyle{definition}
\newtheorem{defn}[thm]{Definition}

\newtheorem{rem}[thm]{Remark}

\numberwithin{equation}{section}

\title[Frame properties of effective stationary sequences]{Frame properties of effective stationary sequences}

\author[Chad Berner]{Chad Berner}
\email{chadberner713@gmail.com}

\author[John E. Herr]{John E. Herr}
\address{Department of Mathematical Sciences, Butler University, Indianapolis, Indiana, USA}
\email{jeherr@butler.edu}

\author[Palle E.T. Jorgensen]{Palle E.T. Jorgensen}
\address{Department of Mathematics, University of Iowa, Iowa City, Iowa, USA}
\email{palle-jorgensen@uiowa.edu}


\date{November 2, 2026}

\begin{document}

\begin{abstract}
The Kaczmarz algorithm is an iterative projection method for reconstructing vectors in a Hilbert space from inner product data. However, sequences that are effective for this algorithm need not be frames and may lack standard frame reconstruction properties. We investigate these properties for  effective stationary sequences, which were classified by Kwapień and Mycielski according to their spectral measures. In particular, we show that the corresponding Kaczmarz reconstruction series always converge unconditionally if and only if the spectral measure is atomic with finitely many point masses. We show that if the spectral measure is atomic, the effective stationary sequence and its auxiliary Parseval frame can be interchanged in the Kaczmarz reconstruction formula. Finally, we construct two effective stationary sequences with singular continuous spectral measures, one of which can be interchanged with its auxiliary Parseval frame in the reconstruction formula, and one of which cannot be interchanged.
\end{abstract}

\maketitle

\section{Introduction}

Consider the problem of recovering vectors in a Hilbert space stably from inner product data. The Kaczmarz algorithm offers an iterative approach to this problem based on successive projections. Originating with Stefan Kaczmarz in 1937 \cite{Kaczmarz1993Approximate} and developed in the infinite-dimensional setting by Kwapień and Mycielski in 2001 \cite{Kwapien2001Kaczmarz}, this method reconstructs vectors using sequences of unit vectors called effective sequences. Frame theory provides another approach to this problem, in which frame inequalities guarantee stable reconstruction from inner product data.

Although these approaches share a common goal, neither effectiveness nor the frame property implies the other. There are frames that are not effective for the Kaczmarz algorithm, and there are effective sequences that are not frames (see \cite{Czaja2013Kaczmarz}, for instance). Thus, the Kaczmarz reconstruction procedure does not by itself determine which frame-implied properties the underlying effective sequence possesses.

This difference motivates a natural question: which properties of frames are shared by effective sequences? In this paper, we investigate this question for the effective stationary sequences characterized by Kwapień and Mycielski. In particular, we show that, apart from trivial cases, the associated reconstruction fails to be unconditional, and the domain of the analysis operator of the effective stationary sequence is trivial. We also investigate whether reconstruction is still guaranteed when the effective stationary sequence and its auxiliary Parseval frame are interchanged, obtaining both positive results and counterexamples.

There is substantial research on the relationship between the Kaczmarz algorithm and frame theory. Kwapień and Mycielski \cite{Kwapien2001Kaczmarz} introduced effective sequences in infinite-dimensional spaces, related effectiveness to the Parseval frame property of an auxiliary sequence, and characterized effective stationary sequences in terms of their spectral measures. Haller and Szwarc \cite{Haller2005Kaczmarz} subsequently characterized effectiveness through a partial isometry condition on a matrix determined by the original sequence.  Szwarc \cite{Szwarc2007Kaczmarz} further developed this connection by showing that every Parseval frame whose first vector has norm one can be realized as the auxiliary sequence of an effective sequence. He also investigated the uniqueness of this correspondence, providing a converse perspective on the relationship between effective sequences and Parseval frames. Czaja and Tanis \cite{Czaja2013Kaczmarz} further investigated the relationship between effectiveness and frame properties, introducing almost effective sequences and studying connections with frames and Riesz bases. Furthermore, Camrud, Camrud, Przybylski, and Weber \cite{Camrud2019Stability} investigated the stability of stationary Kaczmarz reconstructions under additive noise. Their work demonstrates that effectiveness alone does not guarantee stability and establishes stabilization results through a relaxed Kaczmarz algorithm corresponding to Abel summation. These results highlight a further distinction between reconstruction by effective sequences and reconstruction by frames.

In the stationary setting, Herr and Weber \cite{HerWeb2017} applied these ideas to Fourier expansions for singular Borel probability measures, showing that the nonnegative integer exponentials form an effective sequence despite failing the Bessel condition. Additionally, Aboud, Curl, Harding, Vaughan, and Weber \cite{Aboud2020Dual} introduced the dual Kaczmarz algorithm and effective pairs, allowing different sequences to provide the inner product data and the reconstruction directions. Aboud and Weber \cite{Aboud2022Kaczmarz} further developed the Banach space setting, deriving conditions for effectiveness using mixed Gram matrices.

While the focus of our paper is that of frame properties in a general Hilbert space theoretic context, we stress that this in turn is motivated by both the cited classical results involving contexts of Kaczmarz algorithms and recent uses (by the authors and others) of the Kaczmarz algorithm for inverse problems in harmonic analysis. In addition, there is a rich variety of new applications involving the Kaczmarz algorithm as a powerful tool. For this we call the reader’s attention to the following papers: \cite{Yin2024Constrained, Yi2026Recovery, Yuan2025Greedy, LiaoL2024Accelerated, BouRabee2025Accelerated, He2023Dual}. A list of frame terminology and related concepts appearing in this paper is provided at the end.

The aim of this paper is to investigate which frame properties effective stationary sequences possess. Specifically, we study unconditional convergence of the associated reconstruction series, properties of the analysis and synthesis operators, and conditions under which the roles of the effective sequence and its auxiliary Parseval frame can be interchanged in the reconstruction formula.

This paper is organized as follows. In Section 3, we examine the analysis and synthesis operators associated with effective stationary sequences and derive a few consequences for their frame properties. In Section 4, we use techniques from model space theory to show that the associated reconstruction is unconditional only when the spectral measure has finite atomic support. Finally, in Section 5, we establish that the roles of the effective sequence and its auxiliary Parseval frame can be interchanged when the spectral measure is discrete, and we carefully construct an example and a counterexample for singular continuous spectral measures.

\subsection{Note on AI Use}

AI was used extensively in this project for researching, brainstorming, and suggesting proofs. We meticulously scrutinized and rewrote all AI-composed proofs used in this paper, often with substantial editing for the sake of clarity and concision. Among the models most used were Claude Fable 5 and 5.1, Claude Opus 4.8, 5.0, and 5.5, and ChatGPT-5.6 Sol and ChatGPT-6 Sol. Note that updated models were released by Anthropic and OpenAI over the course of this project, contributing to the length of this list.

Along the way, we used as a template the prompt from which ChatGPT-5.6 resolved the Cycle Double Cover Conjecture in July 2026 \cite{OpAI}, though we opted to use five or six parallel subagents rather than 64. However, we believe the explicit prompt call for subagents is becoming obsolete as AI models improve and become more autonomous.

In almost all cases, the AI output was correct but not adequately organized or exposited for human consumption. In a few instances, our own insight was used to finish an argument the AI began but failed to close. Other results are entirely our own.

We have attempted to mark on a result-by-result basis which material was of AI origin. However, the demarcation of lemmas and theorems is itself a product of our own revisions, and the fluidity of AI use makes it impossible to perfectly separate AI output from human output.

In any case, we find the intricacy of the AI-generated proofs to be remarkable given that we used only consumer-grade AI resources.

\section{Preliminaries}

\subsection{Frames and their properties}

We begin by providing the definition of frames:
\begin{defn}
Let $H$ be a Hilbert space and $\{f_{n}\}_{n}\subseteq H$.
\begin{enumerate}
\item If there exist $A,B>0$ such that
\[
A\|f\|^{2}\leq \sum_{n}|\langle f,f_n\rangle|^{2}
\leq B\|f\|^{2}
\]
for all $f\in H$, then
$\{f_n\}_n$ is called a \textbf{frame}\index{Frames}.
\item If
$$\|f\|^{2}=\sum_{n}|\langle f, f_{n}\rangle|^{2}$$
for all $f\in H$, then $\{f_{n}\}_{n}$ is called a \textbf{Parseval frame}\index{Parseval frame property}.
\end{enumerate}
Moreover, if $\{f_{n}\}$ is a frame,
the map $\mathcal{A}:H\to\ell^2$: $f\mapsto \{\langle f, f_{n}\rangle \}_{n}$ is called the \textbf{analysis operator}\index{Analysis and synthesis operators}, its adjoint is called the \textbf{synthesis operator}, and 
we denote by $\mathcal{F}: H\to H$ the positive invertible operator such that
$$\mathcal{F}(f):=\mathcal{A}^\ast\mathcal{A}f=\sum_{n}\langle f, f_{n}\rangle f_{n},$$ which is called the \textbf{frame operator}.
Additionally, $\{\mathcal{F}^{-1}f_{n}\}_{n}$ is a frame such that
$$f=\sum_{n}\langle f, \mathcal{F}^{-1}f_{n}\rangle f_{n}=\sum_{n}\langle f, f_{n}\rangle \mathcal{F}^{-1}f_{n}$$ with convergence in norm for all $f\in H$.

Furthermore, if $\{f_{n}\}$ and $\{g_{n}\}$ are frames such that for all $f\in H$,
$$f=\sum_{n}\langle f, f_{n}\rangle g_{n},$$\index{Stable reconstruction from inner-product data} then $\{f_{n}\}$ and $\{g_{n}\}$ are called \textbf{dual frames}\index{Dual frames}. A frame may have many dual frames, but $\set{\mathcal{F}^{-1}f_n}$ is called the \textbf{canonical dual frame}\index{Canonical dual frame} of $\set{f_n}$. Additionally, because the synthesis and analysis operators associated with both sequences are bounded, by an adjoint calculation, 
$$f=\sum_{n}\langle f,g_{n}\rangle f_{n}$$ also holds with convergence in norm for all $f\in H$. Moreover, this convergence is unconditional.\index{Unconditional and conditional convergence}
\end{defn}

For more information on frames and frame theory, see \cite{christensen2016Introduction}.

\subsection{Inner functions and Model spaces}\index{Inner functions}
We recall the following facts about inner functions and their corresponding measures, which the readers can find in \cite{GarMasRos16} by Garcia, Mashreghi, and Ross.
\begin{defn}
There is a one-to-one correspondence between the inner functions $b$ on $\mathbb{D}$ such that $b(0)=0$ and singular Borel probability measures $\mu$ on $[0,1)$ given by 
\begin{align*}\frac{1+b(z)}{1-b(z)}=\int_{0}^{1}\frac{e^{2\pi ix}+z}{e^{2\pi ix}-z}\,d\mu(x).\end{align*}

Here $\mu$ is called the \textbf{Aleksandrov-Clark measure (with $\alpha=1$)}\index{Aleksandrov-Clark measures} corresponding to the inner function $b$.     
Throughout this paper, we will refer to $\mu$ and $b$ via this correspondence. The integral
\begin{align*}P_{\mu}(z):=&\text{Re}\left(\int_{0}^{1}\frac{e^{2\pi ix}+z}{e^{2\pi ix}-z}\,d\mu(x) \right)\\
=&\int_{0}^{1}\frac{1-\abs{z}^2}{\abs{e^{2\pi ix}-z}^2}\,d\mu(x)\\
=&\frac{1-\abs{b(z)}^2}{\abs{1-b(z)}^2}\end{align*} is the \textbf{Poisson integral of $\mu$}.

Moreover, recall that every inner function factors as $$b(z)=\gamma B(z)S(z),$$ where $\abs{\gamma}=1$, $B(z)$ is a Blaschke product, and $$S(z)=\exp\p{-\int_{0}^{1}\frac{e^{2\pi ix}+z}{e^{2\pi ix}-z}\,d\nu(x)}$$
is the \textbf{singular inner factor of $b(z)$}, where $\nu$ is a finite, positive Borel measure. 
\end{defn}

Note that $\mu$ and $\nu$ are not the same measure. In fact, we have the following:

\begin{prop}
    Let $b(z)$ be an inner function corresponding to $\mu$. Then
    \begin{enumerate}[label=(\roman*)]
        \item $\lim_{r\to 1^{-}}b(re^{2\pi ix})=1$ for $\mu$-almost-every $x\in[0,1)$.
        \item $\lim_{r\to 1^{-}}b(re^{2\pi ix})=0$ for $\nu$-almost-every $x\in [0,1)$. 
         \item $\mu$ is singular with respect to Lebesgue measure.
         \item $\nu$ is singular with respect to Lebesgue measure.
        \item $\mu$ and $\nu$ are  singular with respect to each other.
    \end{enumerate}
\end{prop}
\begin{proof}All five results are well known and can be found in \cite{GarMasRos16} as follows: For $(i)$ and $(iii)$, see Proposition 11.2. For $(ii)$, see the proof of Proposition 7.19. For $(iv)$, see Theorem 2.12. Finally, $(v)$ follows from $(i)$ and $(ii)$ above.
\end{proof}

\begin{defn}
Beurling’s celebrated invariant subspace theorem \newline \cite{Beurling1949Two} implies that the nonzero proper closed subspaces of the Hardy space \(H^2(\mathbb{D})\) invariant under the backward shift are precisely the \textbf{model spaces}, denoted by
\[
\mathcal{H}(b):=H^2(\mathbb{D})\ominus bH^2(\mathbb{D}),
\]
where \(b\) is a nonconstant inner function. With the inner product inherited from \(H^2(\mathbb{D})\), each \(\mathcal{H}(b)\) is a reproducing kernel Hilbert space. Its reproducing kernels, denoted by \(k_z^b\), are given by
\[
k_z^b(w)=\frac{1-\overline{b(z)}b(w)}{1-\overline{z}w},
\qquad z,w\in\mathbb{D},
\]
so that \(f(z)=\langle f,k_z^b\rangle_{\mathcal{H}(b)}\) for every \(f\in \mathcal{H}(b)\).   

The reproducing kernel of the full Hardy space $H^2(\mathbb{D})$ is the Szeg\H{o} Kernel, given by $$k_z(w)=\frac{1}{1-\overline{z}w}.$$
\end{defn}

\subsection{The Kaczmarz algorithm}
\begin{defn}
Let $H$ be a Hilbert space, and let $\{e_n\}_{n=0}^{\infty}\subseteq H$ be a complete sequence of unit vectors in $H$. For an $x\in H$, the \textbf{Kaczmarz algorithm}\index{Kaczmarz algorithm} attempts to recover $x$ by creating a sequence of estimates $\set{x_n}_{n=0}^{\infty}$ recursively as follows:
\begin{align*}x_0=&\ip{x,e_0}e_0\\
x_n=&x_{n-1}+\ip{x-x_{n-1},e_n}e_n.
\end{align*}
In other words, the sequence $\set{x_n}$ consists of successive projections of the previous estimate onto the affine hyperplane containing $x$ perpendicular to $e_n$. If for any $x\in H$, $\lim x_n=x$ in norm, then $\set{e_n}$ is said to be an \textbf{effective sequence}\index{Effective sequences}.

Define the \textbf{auxiliary sequence}\index{Auxiliary Parseval frames} recursively as follows:
\begin{align*}
g_0=&e_0\\
g_n=&e_n-\sum_{k=0}^{n-1}\langle e_n, e_k\rangle g_k \qquad (n\geq 1).\end{align*}
\end{defn}
Then $x_n=\sum_{k=0}^{n}\ip{x,g_k}e_k$. Thus an equivalent formulation of $\set{e_n}$ being effective \index{Effective sequences} is that $$\sum_{n=0}^{\infty}\ip{x,g_n}e_n=x$$ 
for all $x\in H$.

As shown by Kwapie\'{n} and Mycielski in \cite{Kwapien2001Kaczmarz}, $\set{e_n}$ is effective if and only if $\set{g_n}$ is a Parseval frame in $H$. In fact, they proved the following:

\begin{thm*}[Kwapie\'{n} and Mycielski]
If $\{e_n\}$ is a stationary sequence of unit vectors that is linearly dense in $H$, then $\{e_n\}$ is effective if and only if the spectral measure for $\{e_n\}$ is Lebesgue measure or singular. In this case, $\{g_n\}$ is a Parseval frame.
\end{thm*}

By \textbf{stationary}, it is meant that $\ip{e_{n+k},e_{m+k}}_{H}=\ip{e_n,e_m}_{H}$ for all $k\geq0$. By \textbf{spectral measure}\index{Spectral measure}, it is meant the measure $\mu$ on the torus, hereafter identified with $[0,1)$, associated to $\set{e_n}$ by Bochner's Theorem:
\begin{equation}\ip{e_{n+k},e_{m+k}}_{H}=\int_{0}^{1}e^{-2\pi i(m-n)x}\,d\mu(x).
\end{equation}

As observed in \cite{HerWeb2017}, if $H=L^2(\mu)$ for a singular Borel probability measure $\mu$ on $[0,1)$ and $e_n=e^{2\pi inx}$, then $\set{e_n}$ is stationary and complete in $L^2(\mu)$, and $\mu$ is its spectral measure. Thus by the Kwapie\'{n}-Mycielski Theorem, $\set{e_n}$ is effective in $L^2(\mu)$. Any effective sequence with a singular spectral measure $\mu$ is unitarily equivalent to $\set{e^{2\pi inx}}$ in $L^2(\mu)$ via the obvious map. It is the purpose of this paper to study the extent to which $\{e^{2\pi inx}\}_{n=0}^{\infty}$ exhibits frame-like properties in $L^2(\mu)$. We will hereafter take $e_n$ to denote $e^{2\pi inx}$.

\begin{defn}
By the Clark theory \cite{Clark1972One} and \cite{Herr2016Fourier}, the following map, called the \textbf{normalized Cauchy transform}\index{Normalized Cauchy transform}, $V_{\mu}: L^2(\mu)\to \mathcal{H}(b)$ is given by
\begin{equation}\label{NCT_series}\begin{aligned}[V_{\mu}f](z):=&(1-b(z))\int_{0}^{1}\frac{f(x)}{1-ze^{-2\pi ix}}\,d\mu(x)\\
=&\sum_{n=0}^{\infty}\langle f, g_n\rangle z^n,
\end{aligned}\end{equation} and it is unitary.

It is also shown in \cite{Polt93} (and re-demonstrated in\cite{Herr2016Fourier}) that $$[V_\mu f](re^{2\pi i x})\to f(x)$$ as $r\to 1^{-}$ for $\mu$-almost-every $x$. In other words, $f$ is the pointwise radial \textbf{$\mu$ boundary limit of $V_{\mu}f$}\index{$\mu$ boundary of Hardy space functions}.
\end{defn}

Let $P_{\mathcal{H}(b)}$ denote the orthogonal projection of $H^2$ onto $\mathcal{H}(b)$. By equation \eqref{NCT_series}, for any $f\in L^2(\mu)$, 
\begin{align*}\ip{P_{\mathcal{H}(b)}z^n,V_\mu f}_{H^2}=&\ip{z^n,V_\mu f}_{H^2}=\ip{g_n,f}_{L^2(\mu)}=\ip{V_\mu g_n, V_\mu f}_{L^2(\mu)},
\end{align*}
whence it follows that \begin{equation}\label{V_equals_Pgn}V_\mu g_n=P_{\mathcal{H}(b)}z^n.\end{equation}

\begin{defn}
For a singular probability measure $\mu$ on $[0,1)$, we define $T_{\mu}\in B(L^2(\mu))$ by $T_{\mu}=P_{g_0^{\perp}} M_{e^{2\pi i x}}$, where $P_{g_0^{\perp}}$ denotes the orthogonal projection onto  $\operatorname{span}\{g_0\}^{\perp}$.
\end{defn}

We will also use the following, which follows from \cite{Berner2026Operator} and \newline \cite{Christensen2020Frame}:
\begin{thm*}
$T_{\mu}^{n}g_0=g_n$ for all $n$, and for all $f\in L^2(\mu)$,
$T_{\mu}^{n}f\to 0$ as $n\to\infty.$
\end{thm*}

\section{Effective stationary sequences: the straightforward frame properties} We begin by discussing some of the straightforward frame properties of effective stationary sequences. Recall that we may assume we are working with $\{e^{2\pi inx}\}_{n=0}^{\infty}\subseteq L^2(\mu)$ where $\mu$ is a singular probability measure on $[0,1)$.

\subsection{Analysis operator}
\begin{prop}\label{analysisoperator}\index{Analysis and synthesis operators}
The analysis operator corresponding to $\{e^{2\pi inx}\}_{n=0}^{\infty}$ has trivial domain. 
\end{prop}
\begin{proof}
This follows from the F.~and~M.~Riesz theorem; for details, see
the proof of Theorem~7.2 in \cite{Berner2026Continuously}.
\end{proof}

\subsection{Synthesis operator}\index{Analysis and synthesis operators}
Remarkably, as a consequence of the analysis operator, the synthesis operator for $\{e^{2\pi inx}\}_{n=0}^{\infty}$ is, in the sense given below, maximally bad, regardless of the frame-like reproducing property that the $\{e^{2\pi inx}\}_{n=0}^{\infty}$ system possesses:
\begin{cor}
Let $T$ denote the synthesis operator corresponding to $\{e^{2\pi inx}\}_{n=0}^{\infty}$ with domain restricted to $c_{00}$, the space of sequences of finite support. The graph of $T$ is dense in $\ell^2(\mathbb{N})\oplus L^2(\mu)$. 
\end{cor}
\begin{proof}
Let $(\{d_n\},f)\in \ell^2(\mathbb{N})\oplus L^2(\mu)$ and suppose 
$$\langle (\{c_n\}, T\{c_n\}), (\{d_n\},f)\rangle=0 $$ for all $\{c_n\}\in c_{00}$. It follows that
$$\langle T\{c_n\}, f\rangle =\langle \{c_n\}, -\{d_n\}\rangle=\sum_{n}c_n \langle e^{2\pi inx}, f\rangle .$$

Therefore, $\{\langle f, e^{2\pi inx}\rangle\}_{n=0}^{\infty}\in \ell^{2}(\mathbb{N})$, which implies $f=0$ by Proposition \ref{analysisoperator}, and hence $\{d_n\}=0$ as well since $c_{00}$ is dense in $\ell^2(\mathbb{N})$.
\end{proof}

\begin{rem}
Therefore, in the $\ell^2$ sense, the coefficients of a trigonometric  polynomial reveal nothing about the trigonometric  polynomial in $L^2(\mu)$ when $\mu$ is singular.
\end{rem}

However, we know the synthesis operator corresponding to \newline $\{e^{2\pi inx}\}_{n=0}^{\infty}$ can be defined on other spaces, such as the model space it corresponds to by effectiveness. Furthermore, it is shown in \cite{Herr2016Fourier} that there are other closed subspaces of $\ell^2$ where the synthesis operator can be defined. 

Now what if we want the $\{e^{2\pi inx}\}_{n=0}^{\infty}$ series to converge unconditionally? Is it possible to have unconditional convergence of these series without absolute convergence? It turns out that if $\mu$ is not continuous (has atoms), the answer is no.

\begin{prop}
Let $\mu$ be a singular probability measure on $[0,1)$ that is not continuous. If for some $\{c_n\}_{n=0}^{\infty}\subseteq \mathbb{C}$, $\sum_{n=0}^{\infty}c_n e^{2\pi inx}$ converges in $L^2(\mu)$ unconditionally, then the convergence is also absolute.\index{Absolute convergence}
\end{prop}
\begin{proof}
Since $\mu$ is not continuous, it follows that there is a $p\in [0,1)$ with $\mu(\{p\})>0$. Now note that
$$\sum_{n=0}^{\infty}c_n \langle e^{2\pi inx}, \chi_{\{p\}}\rangle=\sum_{n=0}^{\infty}c_n \mu(\{p\})e^{2\pi i np}$$ converges unconditionally and hence absolutely by the Riemann series theorem. Hence $\sum\norm{c_n e^{2\pi inx}}=\sum\abs{c_n}<\infty.$
\end{proof}

However, the answer is yes if the Fourier coefficients of $\mu$ decay fast enough. We know this is possible for example by the following result of Baker and Banaji \cite{Baker2025Polynomial}. We state a weaker version of their result here for our specific use:

\begin{thm*}[Baker and Banaji]
There exist singular probability measures on $[0,1)$ whose Fourier coefficients are in $\ell^{p}$ for some $2<p<\infty$.
\end{thm*}

\begin{lem}[Proof drafted by ChatGPT-5.6, human-revised]\label{analysislem}
Let $\mu$ be a singular probability measure on $[0,1)$ whose Fourier coefficients are in $\ell^{p}$ for some $2<p<\infty$. Then the map $$\mathcal{A}f= \{\langle f, e^{2\pi inx}\rangle\}_{n=0}^{\infty}$$ is bounded from $L^2(\mu)$ into $\ell^{2p}$.
\end{lem}
\begin{proof}
Let $d=\{d_n\}\subseteq \mathbb{C}$ be finitely supported. Then observe that
$$\langle \sum_{n=0}^{\infty}d_n e^{2\pi inx}, \sum_{n=0}^{\infty}d_n e^{2\pi inx}\rangle =\sum_{m,n}d_n \overline{d_m} \hat{\mu}(m-n)\leq$$
$$\sum_{m,n}|d_n| |d_m| |\hat{\mu}(m-n)|\leq \norm{|\hat{\mu}|* |d|}_{2p}\norm{d}_{\frac{2p}{2p-1}}\leq \norm{\hat{\mu}}_{p}\norm{d}_{\frac{2p}{2p-1}}^{2}.$$
Therefore, it follows that the map $\{d_n\}\mapsto \sum_{n=0}^{\infty}d_n e^{2\pi inx}$ is bounded from $\ell^{\frac{2p}{2p-1}}$ to $L^2(\mu)$. In particular, the adjoint map $\mathcal{A}$ where $$\mathcal{A}f= \{\langle f, e^{2\pi inx}\rangle\}_{n=0}^{\infty}$$ is bounded from $L^2(\mu)$ into $\ell^{2p}$.
\end{proof}

\begin{thm}\label{RPex}
There is a singular probability measure $\mu$ on $[0,1)$ and a sequence $\{c_n\}_{n=0}^{\infty}\subseteq \mathbb{C}$ such that $\sum_{n=0}^{\infty}c_n e^{2\pi inx}$ converges unconditionally in $L^2(\mu)$ and not absolutely.\index{Absolute convergence}\index{Unconditional and conditional convergence}
\end{thm}
\begin{proof}
By Lemma \ref{analysislem}, there is a singular probability measure $\mu$ on $[0,1)$ and a $2<p<\infty$ such that the map $$\mathcal{A}f= \{\langle f, e^{2\pi inx}\rangle\}_{n=0}^{\infty}$$ is bounded from $L^2(\mu)$ into $\ell^{2p}$.
Now choose $a\in \left(\frac{2p-1}{2p}, 1 \right]$. Observe that for any unit vector $f\in L^2(\mu)$ and bijection $\sigma: \mathbb{N}\to \mathbb{N} $,
$$\left|\ip{\sum_{n=M}^{N}\sigma(n)^{-a} e^{2\pi i\sigma(n)x},f} \right|\leq \left(\sum_{n=M}^{N}\sigma(n)^{-\frac{2ap}{2p-1}} \right)^{\frac{2p-1}{2p}}\norm{\mathcal{A}}.$$ It follows that $\sum_{n=1}^{\infty}n^{-a} e^{2\pi inx}$ converges unconditionally in $L^2(\mu)$ and not absolutely. 
\end{proof}

We will see soon that if $\mathcal{H}(b)$ (the model space corresponding to $\mu$) is the domain of the synthesis operator corresponding to $\{e^{2\pi inx}\}_{n=0}^{\infty}$, the series do not all converge unconditionally for any $\mu$ for which $L^2(\mu)$ is infinite-dimensional.

\subsection{Relationship between effective stationary sequences and their auxiliary sequences}

Recall that a frame $\{f_n\}$ in a Hilbert space has a canonical dual frame $\{\mathcal{F}^{-1}f_n\}$ \index{Canonical dual frame} where $\mathcal{F}$ is the frame operator. This raises the question for the $\{e^{2\pi inx}\}_{n=0}^{\infty}$ and $\{g_n\}$ system: Can there be even so much as a bounded linear map from $\{e^{2\pi inx}\}_{n=0}^{\infty}$ to $\{g_n\}$ or vise versa? It is easy to see the answer to both questions is no, and this follows easily from Proposition \ref{analysisoperator} and the fact that $\{g_n\}$ is a Parseval frame.

\begin{cor}
There is no bounded linear map from $\{e^{2\pi inx}\}_{n=0}^{\infty}$ to $\{g_n\}$ or vice versa.
\end{cor}

\begin{proof}Suppose there is $L\in B(L^2(\mu))$ such that for all $n$, $Le^{2\pi inx}=g_n$. Then for all $f\in L^2(\mu)$, $$\sum\abs{\ip{L^\ast f,e^{2\pi inx}}}^2=\sum\abs{\ip{f,g_n}}^2=\norm{f}^2<\infty.$$ By Lemma \ref{analysisoperator}, we must have $L^\ast f=0$ for every $f$, which is a contradiction. Conversely, suppose there is $L\in B(L^2(\mu))$ such that for all $n$, $Lg_n=e^{2\pi inx}$. Then for all $f\in L^2(\mu)$, $$\sum\abs{\ip{f,e^{2\pi inx}}}^2=\norm{L^\ast f}^2<\infty,$$ contrary to Lemma \ref{analysisoperator}.
\end{proof}

\section{Unconditional Convergence of the Kaczmarz Series}
Now in this section, we completely answer the question: for which singular probability measures $\mu$ on $[0,1)$ does $\sum_{n=0}^{\infty}\langle f, g_n\rangle e^{2\pi i nx}$ converge unconditionally in $L^2(\mu)$ for all $f\in L^2(\mu)$? We show this can be the case only if $L^2(\mu)$ is finite-dimensional (that is, $\mu$ is a finite sum of point masses).

\subsection{Finite-Dimensional Case}

\begin{thm}\label{finite_theorem}
Let $\mu$ be a probability measure on $[0,1)$ supported on finitely many points. Then for each $f\in L^2(\mu)$, the series
$$\sum_{n=0}^{\infty}\langle f, g_n \rangle e^{2\pi inx}$$ converges absolutely and hence unconditionally.\index{Absolute convergence}\index{Unconditional and conditional convergence}
\end{thm}
\begin{proof}
Utilizing the normalized Cauchy transform and equation \eqref{V_equals_Pgn}, it suffices to show that for all $g(z)\in \mathcal{H}(b)$,
\begin{equation}\label{main}
\sum_{n=0}^{\infty}|\langle g(z), P_{\mathcal{H}(b)}z^n\rangle|=\sum_{n=0}^{\infty}|\langle g(z), z^n\rangle|<\infty.   
\end{equation} Since $\mathcal{H}(b)$ is finite dimensional, it consists of rational functions all with the same poles lying outside of the closed unit disk (see 4.3.3. of \cite{Garcia2023Model} for details). Therefore, there is an $r>1$ such that all $g(z)\in \mathcal{H}(b)$ are analytic on $B(0,r)\subseteq \mathbb{C}$. Choose $1<a<r$. Then for all $g(z)\in \mathcal{H}(b)$,
$$\sum_{n=0}^{\infty}\langle g(z), z^n\rangle a^n$$ converges, so $$|\langle g(z), z^n\rangle| a^n\to 0$$ and equation \eqref{main} follows.
\end{proof}

\subsection{Main Test}

\begin{prop}[Part of a Claude proof draft]\label{WUC_characterization}Let $\mu$ be a singular Borel probability measure on $[0,1)$. Then 

\begin{enumerate}[label=(\roman*)]\item For a fixed $f\in L^2(\mu)$, $\sum\ip{f,g_n}e_n$ converges unconditionally if and only if $\sum\abs{\ip{f,g_n}}\abs{\ip{e_n,h}}<\infty$ for all $h\in L^2(\mu)$.
\item $\sum\ip{f,g_n}e_n$ converges unconditionally for all $f\in L^2(\mu)$ if and only if $\Xi:=\sup\set{\sum_n\abs{\ip{f,g_n}}\abs{\ip{e_n,h}}\,:\,\norm f=\norm h=1}<\infty$.
\end{enumerate}\index{Unconditional and conditional convergence}
\end{prop}

\begin{proof}(i) Fix an $f\in L^2(\mu)$. If $\sum\ip{f,g_n}e_n$ converges unconditionally, then it is weakly unconditionally Cauchy, and the result follows by considering the functionals $\ip{\cdot,h}$. Conversely, if $\sum\abs{\ip{f,g_n}}\abs{\ip{e_n,h}}<\infty$ for all $h\in L^2(\mu)$, then $\sum\ip{f,g_n}e_n$ is weakly unconditionally Cauchy in light of the Riesz Representation Theorem. Then by the \newline Bessaga–Pełczyński $c_0$ Theorem \cite{BesPel58}, $\sum\ip{f,g_n}e_n$ converges unconditionally in $L^2(\mu)$.

(ii) If $\Xi<\infty$, $\sum\ip{f,g_n}e_n$ converges unconditionally for all $f\in L^2(\mu)$ by (i). Conversely, if $\sum\ip{f,g_n}e_n$ converges unconditionally for all $f\in L^2(\mu)$, then by (i), for each $f,h\in L^2(\mu)$, we have \newline $\sum\abs{\ip{f,g_n}}\abs{\ip{e_n,h}}<\infty$, and the result follows by the Uniform Boundedness Principle.
\end{proof}

\begin{defn}Let $\mu$ be a singular probability measure on $[0,1)$. For $0<r\leq\frac{1}{2}$, we define
\begin{align*}B(x,r)&:=\set{y\in[0,1):\min\set{\abs{x-y},1-\abs{x-y}}<r}
\end{align*}
to be the open ball of radius $r$ on the torus and
\begin{align*}
    M(r)&:=\sup_{x\in[0,1)}\set{\mu(B(x,r))}
\end{align*}
to be the supremum of the measures of balls of radius $r$.

Finally, define \begin{align*}G(\alpha):=\sup\set{\sum_{n=0}^{\floor{\alpha}}\abs{\ip{f,g_n}}~:~f\in L^2(\mu),\norm{f}=1}.\end{align*}
\end{defn}

We will also make use of the following result:
\begin{prop}\label{subsequence_prop}\comment{non-AI}Let $\mu$ be a nonzero singular measure on $[0,1)$, and let $c>0$. Let $\set{r_k}$ be a sequence of positive numbers, $0<r_k\leq\frac{1}{2}$, converging to $0$. Then there exist a subsequence $\set{r_{k_j}}$ and a sequence of points $x_j\in[0,1)$ such that $\mu(B(x_j,r_{k_j}))>cr_{k_j}$ for each $j$. Consequently, $M(r_{k_j})>cr_{k_j}$ for each $j$.
\end{prop}

\begin{proof}Suppose the result does not hold. Then there exists some $K$ such that for $k>K$, $\mu(B(x,r_k))\leq cr_k$ for all $x\in [0,1)$. Then since $\mu(\set{x})\leq\mu(B(x,r_k))\leq cr_k$, by letting $k\to\infty$, we see that $\mu$ is atomless. Let $I$ be an interval, and let $\abs{I}$ denote its length. We may cover $I$ with a collection consisting of finitely many disjoint balls $B(x,r_k)$ and finitely many points. Choosing this cover so that no ball is disjoint from $I$, we see that $\mu(I)\leq\frac{c}{2}\abs{I}+2cr_k$. By taking $k\to\infty$, it follows that $\mu(I)\leq\frac{c}{2}\abs{I}$. By a standard argument, this shows that $\mu$ is absolutely continuous, contrary to the assumption. The conclusion follows immediately.
\end{proof}

\begin{lem}[The Ball Test, Claude AI-drafted and human-revised]\label{Ball_Test}
For every $0<r\leq\frac{1}{2}$, $$\Xi\geq\frac{1}{2}\sqrt{M(r)}G\p{\frac{1}{4\pi r}}.$$
\end{lem}

\begin{proof}Fix an $0<r\leq\frac{1}{2}$ and a $B=B(x_0,r)$ with $\mu(B)>0$. Define $h_B:=\frac{1}{\sqrt{\mu(B)}}\chi_{B},$ and note that by construction $h_B$ is a unit vector. 

Let $n$ be an integer with $0\leq n\leq\frac{1}{4\pi r}$, and let $x\in B$. Since the chord between two points on a circle is shorter than the arc between them, we have
\begin{align*}
    \abs{e^{2\pi in(x-x_0)}-1}\leq 2\pi nr\leq \frac{1}{2}.
\end{align*}
Therefore,
\begin{align*}
    \abs{\ip{e_n,h_B}}=&\frac{1}{\sqrt{\mu(B)}}\abs{\int_{B}e^{2\pi inx}\,d\mu}\\
    =&\frac{1}{\sqrt{\mu(B)}}\abs{\int_{B}e^{2\pi inx_0}\,d\mu-\int_{B}e^{2\pi inx_0}-e^{2\pi inx}\,d\mu}\\
    =&\frac{1}{\sqrt{\mu(B)}}\abs{\int_{B}1\,d\mu-\int_{B}1-e^{2\pi in(x-x_0)}\,d\mu}\\
    \geq&\frac{1}{\sqrt{\mu(B)}}\abs{\abs{\int_{B}1\,d\mu}-\abs{\int_{B}1-e^{2\pi in(x-x_0)}\,d\mu}}\\
    =&\frac{1}{\sqrt{\mu(B)}}\abs{\mu(B)-\abs{\int_{B}1-e^{2\pi in(x-x_0)}d\mu}}\\
    \geq&\frac{1}{\sqrt{\mu(B)}}(\mu(B)-\mu(B)(\frac{1}{2}))\\
    =&\frac{1}{2}\sqrt{\mu(B)}.
\end{align*}

Therefore, for $f\in L^2(\mu)$ with $\norm{f}=1$,

\begin{align*}&\sum_{n=0}^{\infty}\abs{\ip{f,g_n}}\abs{\ip{e_n,h_B}}\\
\geq&\sum_{n=0}^{\floor{\frac{1}{4\pi r}}}\abs{\ip{f,g_n}}\abs{\ip{e_n,h_B}}\\
\geq&\sum_{n=0}^{\floor{\frac{1}{4\pi r}}}\abs{\ip{f,g_n}}\frac{1}{2}\sqrt{\mu(B)}.\\
\end{align*}
Consequently,
\begin{align*}\sup_{\substack{\norm{f}=1\\x_0\in[0,1)}}\sum_{n=0}^{\infty}\abs{\ip{f,g_n}}\abs{\ip{e_n,h_{B(x_0,r)}}}\geq\frac{1}{2}\sqrt{M(r)}G\p{\frac{1}{4\pi r}}.
\end{align*}
It follows that \begin{align*}
    \Xi\geq\frac{1}{2}\sqrt{M(r)}G\p{\frac{1}{4\pi r}}.
\end{align*}
\end{proof}

\begin{thm}[The Reduction Theorem, Claude AI-drafted and human-revised]\label{Reduction_Theorem}Let $\mu$ be a singular Borel probability measure on $[0,1)$. If \begin{align*}
    \limsup_{r\to 0}M(r)G\p{\frac{1}{4\pi r}}^2=\infty,
\end{align*}
then there exists an $f\in L^2(\mu)$ for which the series $\sum\ip{f,g_n}e_n$ converges only conditionally.\index{Unconditional and conditional convergence}
\end{thm}

\begin{proof}There exists a decreasing sequence of radii $\set{r_k}$ such that \begin{align*}\lim_{k\to\infty}M(r_k)G\p{\frac{1}{4\pi r_k}}^2=\infty.\end{align*} Consequently,
\begin{align*}
    \lim_{k\to\infty}\frac{1}{2}\sqrt{M(r_k)}G\p{\frac{1}{4\pi r_k}}=\infty,
\end{align*}
and so $\Xi=\infty$ by Lemma \ref{Ball_Test}. The result then follows immediately by Proposition \ref{WUC_characterization}. 
\end{proof}

\subsection{Main Results Obtained from Model Space Theory}

\begin{thm}[Evaluation bound and small modulus points, Claude AI-drafted and human-revised]\comment{Claude AI-written, but improved for simplicity and readability}\label{G_bound}Let $\mu$ be a singular Borel probability measure on $[0,1)$, and let $b$ be its corresponding inner function. Let $w\in\mathbb{D}$, and let $N\geq2$ be a positive integer such that $1-\frac{2}{N}\leq\abs{w}\leq 1-\frac{1}{N}$. If $\abs{b(w)}\leq\frac{\sqrt{6}}{5}$, then 
\begin{align*}\p{G(N)}^2\geq\p{\frac{3}{5}-\frac{\sqrt{6}}{2}\abs{b(w)}}^2 N.
\end{align*}
\end{thm}

\begin{proof} 
Recall that $k_z^b$ denotes the reproducing kernel of $\mathcal{H}(b)$. 

Now,
\begin{align*}
    \norm{k_w^b}=&\sqrt{\ip{k_w^b,k_w^b}}\\
    =&\sqrt{k_w^b(w)}\\
    =&\sqrt{\frac{1-\abs{b(w)}^2}{1-\abs{w}^2}}.
\end{align*}
Therefore,
\begin{align*}
    F:=&\frac{k_w^b}{\norm{k_w^b}}\\
    =&\sqrt{\frac{1-\abs{w}^2}{1-\abs{b(w)}^2}}\cdot\frac{1-\overline{b(w)}b(z)}{1-\overline{w}z}
\end{align*}
is a unit vector in $\mathcal{H}(b)$. Observe that
\begin{align*}
    \hat{F}(n)=&\sqrt{\frac{1-\abs{w}^2}{1-\abs{b(w)}^2}}\p{\overline{w}^n-\overline{b(w)}\widehat{bk_w}(n)}.
\end{align*}
We have that
\begin{align*}
    \sqrt{\sum_{n=0}^{N}\abs{\widehat{bk_w}(n)}^2}\leq&\sqrt{\sum_{n=0}^{\infty}\abs{\widehat{bk_w}(n)}^2}\\
    =&\norm{bk_w}\\
    =&\norm{k_w}\\
    =&\frac{1}{\sqrt{1-\abs{w}^2}}.
\end{align*}
Hence, by the Cauchy-Schwarz Inequality,
\begin{align*}
    &\sum_{n=0}^{N}\abs{\hat{F}(n)}\\
    =&\sqrt{\frac{1-\abs{w}^2}{1-\abs{b(w)}^2}}\sum_{n=0}^{N}\abs{\overline{w}^n-\overline{b(w)}\widehat{bk_w}(n)}\\
    \geq&\sqrt{1-\abs{w}^2}\sum_{n=0}^{N}\p{\abs{w}^n-\abs{b(w)}\abs{\widehat{bk_w}(n)}}\\
    \geq&\sqrt{1-\abs{w}^2}\sum_{n=0}^{N}\abs{w}^n-\sqrt{1-\abs{w}^2}\abs{b(w)}\sum_{n=0}^{N}\abs{\widehat{bk_w}(n)}\\
    \geq&\sqrt{1-\abs{w}^2}\frac{1-\abs{w}^{N+1}}{1-\abs{w}}-\sqrt{N+1}\abs{b(w)}\\
    =&\sqrt{1+\abs{w}}\frac{1-\abs{w}^{N+1}}{\sqrt{1-\abs{w}}}-\sqrt{N+1}\abs{b(w)}\\
    \geq&\frac{1-\abs{w}^{N+1}}{\sqrt{1-\abs{w}}}-\sqrt{N+1}\abs{b(w)}.
\end{align*}

Set $t:=N(1-\abs{w})\in[1,2]$. Since $1-x\leq e^{-x}$ for all real values $x$, we have $\abs{w}\leq e^{-(1-\abs{w})}=e^{-\frac{t}{N}}.$ Therefore,
\begin{align*}
    \abs{w}^{N+1}\leq e^{-t(N+1)/N}\leq e^{-t}.
\end{align*}
Recall that the normalized Cauchy transform $V_\mu$ is a unitary operator from $L^2(\mu)$ to $\mathcal{H}(b)$. Since $N\geq 2$, this shows that
\begin{align*}G(N)\geq&\sum_{n=0}^{N}\abs{\ip{V_\mu^{-1}F,g_n}}\\
=&\sum_{n=0}^{N}\abs{\hat{F}(n)}\\
\geq&\frac{1-\abs{w}^{N+1}}{\sqrt{1-\abs{w}}}-\sqrt{N+1}\abs{b(w)}\\
\geq&\frac{1-e^{-t}}{\sqrt{t/N}}-\sqrt{N+1}\abs{b(w)}\\
\geq&\frac{3}{5}\sqrt{N}-\sqrt{N+1}\abs{b(w)}\\
=&\frac{3}{5}\sqrt{N}-\sqrt{N}\sqrt{1+\frac{1}{N}}\abs{b(w)}\\
\geq&\p{\frac{3}{5}-\frac{\sqrt{6}}{2}\abs{b(w)}}\sqrt{N}.
\end{align*}
Finally, if $\abs{b(w)}\leq\frac{\sqrt{6}}{5}$, then $\p{\frac{3}{5}-\frac{\sqrt{6}}{2}\abs{b(w)}}\geq0$, and so 
\begin{align*}\p{G(N)}^2\geq\p{\frac{3}{5}-\frac{\sqrt{6}}{2}\abs{b(w)}}^2 N,
\end{align*} as we wished to show.
\end{proof}

\begin{rem}The constant $\frac{\sqrt{6}}{5}$ in this proof is not sharp and can be made larger at the cost of simplicity. In particular, the bound $\frac{1-e^{-t}}{\sqrt{t}}\geq\frac{3}{5}$ for $t\in[1,2]$ was used, but in fact $\frac{1-e^{-t}}{\sqrt{t}}\geq\frac{1-e^{-2}}{\sqrt{2}}\approx 0.61141.$ 
\end{rem}

\begin{lem}\label{Small_b_lemma}Let $\mu$ be a singular probability measure on $[0,1)$, and let $b$ be its corresponding inner function. Suppose there exist a $\delta<\frac{\sqrt{6}}{5}$ and a sequence of points $\set{z_k}$ in $\mathbb{D}$ such that $\lim\abs{z_k}=1$ and $\abs{b(z_k)}\leq\delta$ for all $k$. Then there exists an $f\in L^2(\mu)$ such that $\sum\ip{f,g_n}e_n$ converges only conditionally. \index{Unconditional and conditional convergence}
\end{lem}

\begin{proof} 

For each $z_k$, we may select an integer $N_k$ such that $1-\frac{2}{N_k}\leq\abs{z_k}\leq1-\frac{1}{N_k}.$ Indeed, an integer in the interval $\left[\frac{1}{1-\abs{z_k}},\frac{2}{1-\abs{z_k}}\right]$, which is an interval of length no less than 1, suffices. Observe that $\lim_{k\to\infty}N_k=\infty.$ Define $d:=\p{\frac{3}{5}-\frac{\sqrt{6}}{2}\delta}^2>0$. It follows that $(G(N_k))^2\geq dN_k$ by Theorem \ref{G_bound}.

Now select a sequence of positive numbers $\set{r_k}$ in such a way that $\frac{1}{4\pi r_k}=N_k$ for each $k$. Since $N_k\to\infty$, it follows that $r_k\to 0$. By Proposition \ref{subsequence_prop}, for any fixed $c>0$, we may select a subsequence $\set{r_{k_j}}$ such that $M(r_{k_j})>cr_{k_j}$. Therefore,
\begin{align*}M(r_{k_j})\p{G\p{\frac{1}{4\pi r_{k_j}}}}^2\geq cr_{k_j}\frac{d}{4\pi r_{k_j}}=\frac{cd}{4\pi}.
\end{align*}
It follows that 

\begin{align*}\limsup_{r\to 0}M(r)G\p{\frac{1}{4\pi r}}^2\geq\frac{cd}{4\pi}.
\end{align*}
Since $c$ was arbitrary, we have

   \begin{align*}\limsup_{r\to 0}M(r)G\p{\frac{1}{4\pi r}}^2=\infty.
\end{align*}
The result then follows by Theorem \ref{Reduction_Theorem}.
\end{proof}

\begin{rem}Proposition \ref{subsequence_prop} and the observation that \keepmath{\liminf_{r\to0^+}\frac{M(r)}{r} = \infty} for a singular measure, and the consequent \newline Lemma \ref{Small_b_lemma}, were human-provided. The AI version being used at the time failed to close the gap for all singular measures on its own.
\end{rem}

\begin{cor}[Infinitely many zeros]\label{infinite_zeros} Let $\mu$ be a singular probability measure on $[0,1)$, and suppose its corresponding inner function $b$ has infinitely many zeros. Then there exists an $f\in L^2(\mu)$ such that $\sum\ip{f,g_n}e_n$ converges only conditionally.\index{Unconditional and conditional convergence}
\end{cor}

\begin{proof}
     Let $\set{z_k}$ be an enumeration of the zeros of $b$. By the Blaschke condition, $\lim\abs{z_k}=1$. The result follows by applying Lemma \ref{Small_b_lemma}.
\end{proof}

\begin{cor}[Nonconstant $S(z)$, Claude AI-drafted and human-revised]\label{nonconstant_S}Let $\mu$ be a singular probability measure on $[0,1)$, and suppose its corresponding inner function $b(z)=\gamma B(z)S(z)$ has singular inner factor $S(z)\not\equiv 1$. Then there exists an $f\in L^2(\mu)$ such that $\sum\ip{f,g_n}e_n$ converges only conditionally.\index{Unconditional and conditional convergence}
\end{cor}

\begin{proof}
    Write \begin{align*}S(z)=\exp\p{-\int_{0}^{1}\frac{e^{2\pi ix}+z}{e^{2\pi ix}-z}\,d\nu(x)}.
    \end{align*}
    Since $S(z)\not\equiv1$, $\nu\neq0$. We also know that $\nu$ is singular with respect to Lebesgue measure.

    For an $x\in[0,1)$ and $0<h\leq\frac{1}{2}$, set $w=(1-h)e^{2\pi ix}$. If $y\in B(x,h)$, then we have
    \begin{align*}
        \abs{e^{2\pi iy}-w}\leq&\abs{e^{2\pi iy}-e^{2\pi ix}}+\abs{e^{2\pi ix}-w}\\
        \leq&2\pi h+h\\
        =&(2\pi+1)h.
    \end{align*}
We observe that $1-\abs{w}^2=1-(1-h)^2=h(2-h)\geq h$. Therefore, the Poisson integral of $\nu$ satisfies
\begin{align*}P_\nu(w)=&\int_{0}^{1}\frac{1-\abs{w}^2}{\abs{e^{2\pi it}-w}^2}\,d\nu(t)\\
\geq&\frac{h}{(2\pi+1)^2h^2}\nu(B(x,h))\\
=&\frac{1}{(2\pi+1)^2}\frac{\nu(B(x,h))}{h}.
\end{align*}
Now let $\set{h_k}$ be a sequence such that $0<h_k\leq\frac{1}{2}$ and $\lim h_k=0$. By Proposition \ref{subsequence_prop}, there exist a subsequence $\set{h_{k_j}}$ and points $x_{j}\in[0,1)$ such that \begin{align*}
    \frac{\nu(B(x_j,h_{k_j}))}{h_{k_j}}>(2\pi+1)^2\ln(3).
\end{align*}
Consequently, for $w_j:=(1-h_{k_j})e^{2\pi ix_j}$, we have \begin{align*}P_\nu(w_j)\geq\ln(3),
\end{align*} and therefore \begin{align*}
    \abs{b(w_j)}\leq&\abs{S(w_j)}\\
    =&e^{-P_\nu(w_j)}\\
    \leq &e^{-\ln(3)}\\
    =&\frac{1}{3}<\frac{\sqrt{6}}{5}.
\end{align*}
   The proof is finished by applying Lemma \ref{Small_b_lemma}.
\end{proof}

\begin{thm}[Main Theorem]\label{Mainthm} Let $\mu$ be a singular probability measure on $[0,1)$. Then the series $\sum_{n}\ip{f,g_n}e^{2\pi i nx}$ converges unconditionally for all $f\in L^2(\mu)$ if and only if $\mu$ is a finite sum of point masses.\index{Unconditional and conditional convergence}
\end{thm}

\begin{proof}If $\mu$ is a finite sum of point masses, every sum $\sum\ip{f,g_n}e^{2\pi inx}$ converges absolutely and hence unconditionally by Theorem \ref{finite_theorem}.

So suppose $\mu$ is not a finitely supported atomic measure, and let $b=\gamma B(z)S(z)$ be its corresponding inner function. It is known that $b$ is a pure finite Blaschke product if and only if $\mu$ is a finitely supported atomic measure, see \cite{Garcia2023Model} for example. Therefore, either $B(z)$ is an infinite Blaschke product or $S(z)\not\equiv 1$. Applying Corollary \ref{infinite_zeros} in the former case and Corollary \ref{nonconstant_S} in the latter, we see that for some $f\in L^2(\mu)$, $\sum\ip{f,g_n}e^{2\pi inx}$ converges conditionally.
\end{proof}

\section{Adjoint Fourier Frame Property for Effective Stationary Sequences}\index{Adjoint frame property for effective stationary sequences}
When $\mu$ is a singular Borel probability measure, we know it is always the case that for $f\in L^2(\mu)$, $$f=\sum_{n=0}^{\infty}\ip{f,g_n}e^{2\pi inx},$$ so naturally, we consider the problem of which singular probability measures on $[0,1)$ have 
\begin{equation}\label{mainadjointeq}
f=\sum_{n=0}^{\infty}\langle f, e^{2\pi inx}\rangle g_n   
\end{equation}
 in $L^2(\mu)$ for all $f\in L^2(\mu)$, motivated by the fact that every pair of dual frames has this property.

In this section, we give a partial answer to this problem, demonstrating that \eqref{mainadjointeq} always holds for discrete measures, and that there exist singular continuous measures for which it holds and others for which it does not hold. However, first we reduce the problem to convergence of one series.

\begin{thm}\label{adjointsimplify}
Let $\mu$ be a singular probability measure on $[0,1)$. Then equation \eqref{mainadjointeq} holds if and only if $$\sum_{n=1}^{\infty}\langle 1, e^{2\pi inx}\rangle g_n=0.$$ \index{Adjoint frame property for effective stationary sequences}
\end{thm}
\begin{proof}
We assume that $$\sum_{n=1}^{\infty}\langle 1, e^{2\pi inx}\rangle g_n=0$$ since the other direction is easy to check.

First, we prove by induction that $$\sum_{n=0}^{\infty}\langle e^{2\pi ikx}, e^{2\pi inx}\rangle g_n=e^{2\pi ikx}$$ for all $k\geq 0.$

The base case easily follows by assumption.

Now suppose that for some $k\geq 0$, $$\sum_{n=0}^{\infty}\langle e^{2\pi ikx}, e^{2\pi inx}\rangle g_n=e^{2\pi ikx}.$$ Then consider:

$$\sum_{n=0}^{N}\langle e^{2\pi i(k+1)x}, e^{2\pi inx}\rangle g_n=\langle e^{2\pi i(k+1)x},1\rangle g_0+ \sum_{n=0}^{N-1}\langle e^{2\pi ikx}, e^{2\pi inx}\rangle g_{n+1}\to$$
$$\langle e^{2\pi i(k+1)x},1\rangle g_0+T_{\mu}e^{2\pi ikx}=e^{2\pi i(k+1)x}.$$ 

Now note that the linear maps $R_{N}: L^2(\mu)\to L^2(\mu)$ given by $$R_{N}(f)=\sum_{n=0}^{N}\langle f, e^{2\pi inx}\rangle g_n$$ are the adjoints of the maps $f\mapsto \sum_{n=0}^{N}\langle f, g_n\rangle e^{2\pi inx}$ given by the Kaczmarz algorithm, which are uniformly bounded by the uniform boundedness principle. Therefore, pointwise convergence on a complete sequence implies pointwise convergence on $L^2(\mu)$.
\end{proof}

\begin{lem}\label{adjointsimplifylem}
Let $\mu$ be a singular probability measure on $[0,1)$. Then
for all $N\geq 1$,
$$T_{\mu}^{N+1}e^{-2\pi i (N+1)x}=-\sum_{n=1}^{N}\langle 1, e^{2\pi i nx}\rangle g_{n}.$$
\end{lem}
\begin{proof}
We start with the $N=1$ case. Observe:
$$T_{\mu}^2(e^{-4\pi ix})=P_{g_0^{\perp}}M_{e^{2\pi ix}}P_{g_0^{\perp}}e^{-2\pi ix}=P_{g_0^{\perp}}(g_0-\langle e^{-2\pi ix}, g_0\rangle e^{2\pi ix})=$$
$$-\langle 1, e^{2\pi ix}\rangle T_{\mu}g_0.$$

Now suppose the equation holds for some $N\geq 1.$ Observe:
$$T_{\mu}^{N+1}e^{-2\pi i(N+1)x}-\langle 1, e^{2\pi i (N+1)x}\rangle g_{N+1}=$$
$$T_{\mu}^{N+1}e^{-2\pi i(N+1)x}-\langle e^{-2\pi i (N+1)x}, g_0\rangle T_{\mu}^{N+1}g_0=T_{\mu}^{N+1}(T_{\mu}e^{-2\pi i (N+2)x}).$$
The result now follows by assuming the equation holds for $N$.
\end{proof}

\begin{defn}\label{adjointsimplifycor}Let $\mu$ be a singular probability measure on $[0,1)$. Define $$A_N(\mu):=\norm{T_{\mu}^{N+1}e_{-(N+1)}}_{L^2(\mu)}.$$
Then in light of Lemma \ref{adjointsimplifylem}, $\mu$ has the \textbf{adjoint Fourier frame property}, that is, equation \eqref{mainadjointeq} holds, if and only if $$\lim_{N\to\infty}A_N(\mu)=0.$$
\end{defn}

\begin{thm}[Chatgpt 5.6 conclusion and proof idea, human edited]\label{discretethm}
Let $\mu$ be a discrete probability measure on $[0,1)$. Then $\mu$ has the adjoint Fourier frame property.\index{Adjoint frame property for effective stationary sequences}
\end{thm}
\begin{proof}
Recall that $\|T_{\mu}\|\leq 1$ and for all $f\in L^2(\mu)$, $\norm{T_{\mu}^{n}f}\to 0.$ 

Now we claim that $\{e^{-2\pi i(N+1)x}\}_{N=1}^{\infty}\subseteq L^2(\mu)$ is pre-compact. To see this, let $\{\lambda_{j}\}_{j=0}^{\infty}$ be a countable set on which $\mu$ is carried and consider any subsequence $\{e^{-2\pi i(N_{k}+1)x}\}$. By a classic diagonalization argument, there is a further subsequence $\{e^{-2\pi i(N_{k_{m}}+1)\lambda_{j}}\}$ that converges pointwise for every $j$. Additionally, $\{e^{-2\pi i(N_{k_{m}}+1)x}\}$ converges in $L^2(\mu)$ by dominated convergence since each function is bounded by 1.

Now, using Definition \ref{adjointsimplifycor}, suppose for the sake of contradiction  that there is an $\epsilon>0$ and a subsequence $N_{k}$ such that $$A_{N_k}(\mu)=\norm{T_{\mu}^{N_k+1}e^{-2\pi i(N_{k}+1)x}}\geq\epsilon$$ for all $k$. We may assume by taking a further subsequence that $$e^{-2\pi i(N_{k}+1)x}\to g$$ for some $g\in L^2(\mu)$. Finally, observe
$$\norm{T_{\mu}^{N_k+1}e^{-2\pi i(N_{k}+1)x}}\leq $$
$$\norm{T_{\mu}^{N_k+1}(e^{-2\pi i(N_{k}+1)x}-g)}+\norm{T_{\mu}^{N_k+1}g},$$ which will be smaller than $\epsilon$ for $k$ large enough, and this is a contradiction. Thus, $\lim A_N(\mu)=0$, and so equation \eqref{mainadjointeq} holds.

\end{proof}
\begin{rem}
It can be shown that if $\{e^{-2\pi i(N+1)x}\}_{N=1}^{\infty}\subseteq L^2(\mu)$ is pre-compact then $\mu$ is discrete by using Wiener's Lemma and noting that if $\{e^{-2\pi i(N+1)x}\}_{N=1}^{\infty}\subseteq L^2(\mu)$ is pre-compact, then $\{e^{-2\pi i(N+1)x}\}_{N=1}^{\infty}\subseteq L^2(\mu_{c})$ is pre-compact where $\mu_{c}$ denotes the continuous part of $\mu$.
\end{rem}

Unsurprisingly, we show that the adjoint frame property is invariant under rotation of the measure:

\begin{thm}
Let $\mu$ be a singular probability measure on $[0,1)$. Fix $\theta\in (0,1)$, and define 
$\mu_{\theta}$ to be shift of $\mu$:
$$\int_{0}^{1}fd\mu_{\theta}=\int_{0}^{1}f(x+\theta \ \operatorname{mod} \ 1) d\mu.$$ Then $\mu$ has the adjoint Fourier frame property if and only if $\mu_{\theta}$ does.\index{Adjoint frame property for effective stationary sequences}
\end{thm}
\begin{proof}
Fix $\theta$ and let $\{g_{n}^{\theta}\}_{n=0}^{\infty}$ denote the auxiliary sequence corresponding to $\mu_{\theta}$. We first prove by induction that for all $n\geq 0$,
$$g_{n}^{\theta}(x)=e^{2\pi in \theta}g_n(x-\theta \ \operatorname{mod} \ 1).$$

The base case is obvious.

For the inductive step, observe that for $n\geq 1$,
\begin{align*}&g_{n}^{\theta}\\
=&
e^{2\pi i n\theta}e^{2\pi i n(x-\theta)}-\sum_{j=0}^{n-1}e^{2\pi i(n-j)\theta}\langle e^{2\pi i nx}, e^{2\pi i jx}\rangle_{\mu}e^{2\pi i j\theta}g_{j}(x-\theta \ \operatorname{mod} \ 1)\\
=&e^{2\pi in \theta}g_n(x-\theta \ \operatorname{mod} \ 1).\end{align*}

Therefore, for all $N\geq 1$,
$$\sum_{n=1}^{N}\hat{\mu}_{\theta}(n)g_{n}^{\theta}=\sum_{n=1}^{N}\hat{\mu}(n)g_n(x-\theta \ \operatorname{mod} \ 1), $$ and the result easily follows.
\end{proof}

\input{SingularContinuousExamples} 

\printindex
\printbibliography

\end{document}

%% file: SingularContinuousExamples.tex
Throughout the rest of the paper, we will use 
$$e_n:=e^{2\pi inx}$$ for shorthand notation.
\begin{lem}[Finitely atomic measures, Claude Opus 5.5 AI-drafted and human-revised]\label{lem:finite}
Let $\nu=\sum_{i=1}^mw_i\delta_{x_i}$ with distinct $x_i\in[0,1)$ and weights $w_i>0$, with $\sum w_i=1$. Then $T_\nu$
has no eigenvalue of modulus $1$. Consequently $\norm{T_\nu^n}$ is nonincreasing in $n$ and tends
to $0$,
\[
  \tau(\nu):=\sum_{n=0}^\infty\norm{T_\nu^n}\ <\ \infty ,
\]
and $A_N(\nu)\le\norm{T_\nu^{N+1}}$ for every $N\ge0$.
\end{lem}

\begin{proof}
Suppose $T_\nu f=\lambda f$ with $\abs\lambda=1$ and $f\ne0$. Since $e_1f=T_\nu f+P_0(e_1f)$
orthogonally and $\norm{e_1f}=\norm f=\norm{T_\nu f}$, we get $P_0(e_1f)=0$. Hence
$e_1f=T_\nu f=\lambda f$, so $(e^{2\pi ix_i}-\lambda)f(x_i)=0$ for every $i$. The numbers
$e^{2\pi ix_i}$ are distinct, so $f$ vanishes at every atom except possibly one atom $x_{i_0}$, where
$e^{2\pi ix_{i_0}}=\lambda$. But then $P_0(e_1f)=\bigl(\int e_1f\,d\nu\bigr)1=\lambda f(x_{i_0})w_{i_0}\cdot1\ne0$,
a contradiction.

$L^2(\nu)$ has dimension $m$, so the spectral radius of $T_\nu$ is the largest modulus of an
eigenvalue, which is less than $1$. By Gelfand's formula, $\norm{T_\nu^n}^{1/n}$ converges to this spectral
radius, so $\sum_n\norm{T_\nu^n}<\infty$. Monotonicity follows from
$\norm{T_\nu^{n+1}}\le\norm{T_\nu}\,\norm{T_\nu^n}\le\norm{T_\nu^n}$. Finally
$$A_N(\nu)\le\norm{T_\nu^{N+1}}\,\norm{e_{-(N+1)}}=\norm{T_\nu^{N+1}}.$$
\end{proof}

\subsection{A Singular Continuous Measure With the Adjoint Fourier Frame Property}\index{Adjoint frame property for stationary effective sequences}

The reader may be familiar with the ternary Cantor measure construction. In this construction, one starts with the unit interval $[0,1]$. Then from it, the open middle-third interval is removed, and from each of the remaining intervals the open middle third is removed, and so on, until in the limit all that remains is the ternary or ``middle-thirds'' Cantor set. The ternary Cantor measure is then the $\frac{\ln(2)}{\ln(3)}$-dimensional Hausdorff measure restricted to the Cantor set. Equivalently, the Cantor measure can be thought of as the pushforward of the $\frac{1}{2}$-$\frac{1}{2}$ weighted Bernoulli measure on the space of sequences $\set{0,1}^\mathbb{N}$ under the map $c\mapsto\sum_{n=1}^{\infty}\frac{2c_n}{3^n}$. Or, one can think of it as the $\frac{1}{2}$-$0$-$\frac{1}{2}$ weighted Bernoulli measure on the space $\set{0,1,2}^{\mathbb{N}}$ of ternary digit sequences, pushedforward to the corresponding number in $[0,1]$ by $c\mapsto\sum_{n=1}^{\infty}\frac{c_n}{3^n}$.

If at the $kth$ step of the Cantor measure construction an atom is placed at, say, the left endpoint of each of the $2^k$ remaining subintervals, and an atomic measure $\nu_k$ created by equally weighting these atoms, the measures $\nu_k$ will converge weakly to the Cantor measure.

The two constructions below follow in a similar vein, but with two notable differences from the Cantor measure construction: First, the scaling factors change at each step. Whereas in the Cantor measure construction one always divides the remaining intervals into thirds, in the following constructions the subdivision is variable at each step, each remaining interval of length $1/Q_j$ being broken into $Q_{j+1}/{Q_j}$ pieces. Second, the probability weights of the first construction are also variable at each step.

We begin by creating a singular continuous measure that has the adjoint Fourier frame property.

\begin{defn}[Claude Opus 5.5 AI-drafted and human-revised]\label{with_adjoint_property_measure}

Let $\set{p_k}$ be a sequence of numbers such that $p_k\in(0,\frac{1}{2}]$, $\lim p_k=0$, and $\sum p_k=\infty$. We define a sequence $\set{Q_k}$ of positive integers and a sequence $\set{\nu_k}$ of finitely supported atomic measures on $[0,1)$ as follows:

Let $Q_0:=1$, and let $\nu_0=\delta_0$.

Once $Q_0,\ldots, Q_{k}$ and $\nu_0,\ldots,\nu_k$ have been defined, let $Q_{k+1}$ be a multiple of $Q_k$ large enough to satisfy the following four properties:

\begin{enumerate}[label=(\roman*)]\item $Q_{k+1}\geq 4 Q_k$;
\item $Q_{k+1}\geq 2\pi \tau(\nu_{k})$;
\item $\norm{T_{\nu_{k}}^{Q_{k+1}}}\leq\sqrt{p_{k+1}}$;
\item $Q_{k+1}\geq(k+1)Q_k(Q_k+1)$.
\end{enumerate}

Note that property (iii) is possible due to Lemma \ref{lem:finite}.

Then define $$\nu_{k+1}:=\sum_{c\in\set{0,1}^{k+1}}\p{\prod_{j=1}^{k+1}\p{c_jp_j+(1-c_j)(1-p_j)}}\delta_{\sum_{j=1}^{k+1}\frac{c_j}{Q_j}}.$$

Observe that $\nu_{k+1}$ is an atomic measure supported on $2^{k+1}$ points. The sequence $\set{\nu_k}$ converges weakly to a measure $\mu$ on $[0,1)$. 

To be more explicit about this, let $$\Omega=\prod_{j=1}^{\infty}\set{0,1,\ldots,\frac{Q_j}{Q_{j-1}}-1}$$ denote the space of sequences of digits used in the above construction, where we remove from $\Omega$ any sequence with a tail of the form $\frac{Q_j}{Q_{j-1}}-1$ in order to make the digits representing a number unique. We may define a measure $m$ on $\Omega$ by specifying it on cylinder sets as follows: Fix a $k$ and a $c\in\prod_{j=1}^{k}\set{0,1,\ldots,\frac{Q_{j}}{Q_{j-1}}-1}$, and set
$$m\p{\set{\omega\in\Omega:\omega_n=c_n, n\leq k}}=\prod_{j=1}^{k}\begin{cases}(1-p_j)&c_j=0\\
p_j&c_j=1\\
0&2\leq c_j
\end{cases}.$$
Then extend $m$ to a measure on $\Omega$ by the Kolmogorov Extension Theorem.

Let $\Phi:\Omega\to[0,1)$ be defined by $\Phi 
(\omega)=\sum_{j=1}^{\infty}\frac{\omega_j}{Q_j}.$
Let $\Psi_k:\Omega\to[0,1)$ be defined by $\Psi_k (\omega)=\sum_{j=1}^{k}\frac{\omega_j}{Q_j}.$ Note that $\Phi$ maps the digit sequence representation of a number to that number, whereas $\Psi_k$ truncates a digit representation to the first $k$ digits and then maps it. The image $\Phi(\Omega)$ is $[0,1)$, whereas the image $\Psi_k(\Omega)$ consists of a finite set of points. The interval $I_c$ is defined as the image of the cylinder set mentioned earlier: $I_c=\Phi\p{\set{\omega\in\Omega:\omega_n=c_n, n\leq k}}=\left[\sum_{j=1}^{k}\frac{c_j}{Q_j},\sum_{j=1}^{k}\frac{c_j}{Q_j}+\frac{1}{Q_k}\right)$.

We may define a measure $\mu$ on $[0,1)$ by
$$\mu(B):=m(\Phi^{-1}(B))$$ for any Borel set $B\subseteq[0,1),$
where $\Phi^{-1}(B)$ denotes the preimage of $B$. Likewise, observe that the measures $\nu_k$ satisify
$$\nu_{k}(B)=m(\Psi_k^{-1}(B)).$$ These observations will be useful to us later.

It is not hard to see that $\mu$ is the measure to which $\set{\nu_k}$ converges weakly, and indeed this is all a typical Bernoulli pushforward construction (see \cite{BHR19}, \cite{Kak48}, and a series of lectures given by Boris Solomyak at Bar-Ilan in 2018). 

\end{defn}

\begin{lem}[Claude Opus 5.5 AI-drafted and human-revised]\label{singular_continuous}
    The measure $\mu$ of Definition \ref{with_adjoint_property_measure} is a singular continuous probability measure.
\end{lem}
\begin{proof}That $\mu$ is a probability measure follows from the clear fact that $m$ is a probability measure.

Fix a $k\geq0$. Observe that $\mu$ is carried on the $2^k$ intervals $I_c$ with $c\in\set{0,1}^k$. Each such interval has length $\frac{1}{Q_k}$, so $\mu$ is carried on a set of Lebesgue measure $\frac{2^k}{Q_k}\leq\frac{2^k}{4^k}=\frac{1}{2^k}$ by property $(i)$. Letting $k\to\infty$, we see that $\mu$ is singular.

Now, $p_j\leq1-p_j$, and so it follows that for any\\$c\in\prod_{j=1}^{k}\set{0,1,\ldots,\frac{Q_j}{Q_{j-1}}-1},$
\begin{align*}\mu(I_c)\leq\prod_{j=1}^{k}(1-p_j).
\end{align*}

Setting $S_k=\prod_{j=1}^{k}(1-p_j)$, by the fact that $\ln(x)\leq x-1$, we have 
\begin{align*}\ln(S_k)=&\sum_{j=1}^{k}\ln(1-p_j)\\
\leq&\sum_{j=1}^{k}\p{(1-p_j)-1}\\
=&-\sum_{j=1}^{k}p_j.
\end{align*}
This last expression diverges to $-\infty$ by assumption, so $S_k\to 0$. Since for any $k$ an atom must lie in one of the intervals $I_c$, we see that in fact $\mu$ has no atoms.
\end{proof}

Despite the scaling factors changing at each step, there is a recursive relationship between $\nu_{k-1}$ and $\nu_k$, the latter simply consisting of all the atoms of the former plus shifted copies, weighted according to $1-p_k$ and $p_k$. Therefore, $L^2(\nu_k)$ is unitarily equivalent to $L^2(\nu_{k-1})\oplus L^2(\nu_{k-1})$, and we now construct a unitary map $F_k$ that will take advantage of this. This will be useful later in placing bounds on the difference between consecutive $A_N(\nu_k)$.

\begin{lem}[Claude Opus 5.5 AI-drafted and human-revised]\label{Fk_Operator}The map $F_k:L^2(\nu_{k})\to L^2(\nu_{k-1})\oplus L^2(\nu_{k-1})$ defined by 
\begin{align*}F_kf\\
=&\p{(1-p_{k})f(x)+p_{k}f\p{x+\frac{1}{Q_{k}}},\sqrt{p_k(1-p_k)}\p{f\p{x+\frac{1}{Q_k}}-f(x)}}\\
=&:(u,v)\end{align*}
is a unitary operator.
\end{lem}

\begin{proof}Linearity of $F_k$ is obvious. Observe by the definition of $\nu_k$ that 
\begin{equation}\label{nu_recursion_equation}\begin{aligned}\int f(x)\,d\nu_{k}(x)=(1-p_{k})&\int f(x)\,d\nu_{k-1}(x)\\
+p_{k}&\int f(x+\frac{1}{Q_{k}})\,d\nu_{k-1}(x).
\end{aligned}
\end{equation}

Set \begin{align*}u(x)=&(1-p_{k})f(x)+p_{k}f\p{x+\frac{1}{Q_{k}}},\\
v(x)=&\sqrt{p_k(1-p_k)}\p{f\p{x+\frac{1}{Q_k}}-f(x)},
\end{align*}
so that $F_kf=(u,v)$. We have $$-\frac{\sqrt{p_k}}{\sqrt{1-p_k}}v(x)=-p_k\p{f\p{x+\frac{1}{Q_k}}-f(x)}.$$
Adding this to $u(x)$, we have
$$f(x)=u(x)-\frac{\sqrt{p_k}}{\sqrt{1-p_k}}v(x).$$
By a similar computation,
$$f\p{x+\frac{1}{Q_k}}=u(x)+\frac{\sqrt{1-p_k}}{\sqrt{p_k}}v(x).$$
Then
\begin{align*}&(1-p_k)\abs{f(x)}^2\\
=&(1-p_k)(\abs{u(x)}^2-\frac{\sqrt{p_k}}{\sqrt{1-p_k}}2\text{Re}\p{u(x)\overline{v(x)}}+\frac{p_k}{1-p_k}\abs{v(x)}^2)
\end{align*}
and
\begin{align*}&p_k\abs{f\p{x+\frac{1}{Q_k}}}^2\\
=&p_k(\abs{u(x)}^2+\frac{\sqrt{1-p_k}}{\sqrt{p_k}}2\text{Re}\p{u(x)\overline{v(x)}}+\frac{1-p_k}{p_k}\abs{v(x)}^2).
\end{align*}
Adding these together, we get 
\begin{align*}(1-p_k)\abs{f(x)}^2+p_k\abs{f\p{x+\frac{1}{Q_k}}}^2=\abs{u(x)}^2+\abs{v(x)}^2.
\end{align*}
So by Equation (\ref{nu_recursion_equation}), it follows that $F_k$ is an isometry. Since the dimension of both its domain and codomain is $2^k$, it is therefore unitary.
\end{proof}

\begin{rem}Observe that $F_k1=(1,0)$, and $\int f(x)\,d\nu_k=\int u(x)\,d\nu_{k-1}(x),$ where $F_kf=(u,v).$
\end{rem}

\begin{lem}[Claude Opus 5.5 AI-drafted and human-revised]\label{FM_lemma}Let $\omega=e^{2\pi i(1/Q_k)}$. If 
$$F_kf=(u,v),$$
then 
$$F_kM_{e_1}f=\p{M_{e_1}(au+sv),M_{e_1}(su+dv)},$$
where 
\begin{align*}a=&1-p_k+p_k\omega\\
s=&\sqrt{p_k(1-p_k)}(\omega-1)\\
d=&p_k+(1-p_k)\omega.
\end{align*}
\end{lem}

\begin{proof}Let $F_kf=(u,v)$. The first coordinate of $F_k(e_1f)$ is by definition
\begin{align*}&(1-p_k)e_1f+p_ke^{2\pi i\p{x+1/Q_k}}f(x+1/Q_k)\\
=&(1-p_k)e_1f+p_k\omega e_1f(x+1/Q_k)\\
=&e_1\p{(1-p_k)f+p_k\omega f(x+1/Q_k)}.
\end{align*}
So for the first coordinate, it remains only to show that $au+sv=(1-p_k)f+p_k\omega f(x+1/Q_k)$. Observe that
\begin{align*}&au+sv\\
=&(1-p_k+p_k\omega)u+\sqrt{p_k(1-p_k)}(\omega-1)v\\
=&(1-p_k+p_k\omega)((1-p_k)f(x)+p_kf(x+1/Q_k))\\
&+\sqrt{p_k(1-p_k)}(\omega-1)\sqrt{p_k(1-p_k)}(f(x+1/Q_k)-f(x))\\
=&(1-p_k+p_k\omega)((1-p_k)f(x)+p_kf(x+1/Q_k))\\
&+p_k(1-p_k)(\omega-1)(f(x+1/Q_k)-f(x))\\
=&(1-p_k+p_k\omega-p_k+p_k^2-p_k^2\omega)f(x)\\
&+(p_k-p_k^2+p_k^2\omega)f(x+1/Q_k)\\
&+(p_k\omega-p_k-p_k^2\omega+p_k^2)f(x+1/Q_k)\\
&-(p_k\omega-p_k-p_k^2\omega+p_k^2)f(x)\\
=&(1-p_k)f(x)+p_k\omega f(x+1/Q_k).
\end{align*}
So the result holds for the first coordinate. The second coordinate of $F_k(e_1f)$ is by definition

\begin{align*}
    &\sqrt{p_k(1-p_k)}\p{e^{2\pi i(x+1/Q_k)}f(x+1/Q_k)-e_1 f}\\
    =&\sqrt{p_k(1-p_k)}(\omega e_1 f(x+1/Q_k)-e_1 f)\\
    =&e_1\sqrt{p_k(1-p_k)}(\omega f(x+1/Q_k)-f(x)).
\end{align*}
Thus it remains only to show that\\ $su+dv=\sqrt{p_k(1-p_k)}(\omega f(x+1/Q_k)-f(x)).$ Observe that
\begin{align*}&su+dv\\
=&\sqrt{p_k(1-p_k)}(\omega-1)\p{(1-p_k)f(x)+p_kf(x+1/Q_k)}\\
&+(p_k+(1-p_k)\omega)\sqrt{p_k(1-p_k)}(f(x+1/Q_k)-f(x))\\
=&\sqrt{p_k(1-p_k)}[(\omega-\omega p_k-1+p_k-p_k-(1-p_k)\omega) f(x)\\
&+(\omega p_k-p_k+p_k+(1-p_k)\omega)f(x+1/Q_k)]\\
=&\sqrt{p_k(1-p_k)}[-f(x)+\omega f(x+1/Q_k)].
\end{align*}
This concludes the proof.
\end{proof}

\begin{lem}[human-provided clarification]\label{tilde_T}Define $\tilde{T}_{\nu_k}:=F_k T_{\nu_k}F_{k}^{-1}$ to be the action of $T_{\nu_k}$ on the space $L^2(\nu_{k-1})\oplus L^2(\nu_{k-1})$. Then $\tilde{T}_{\nu_k}$ is a contraction by unitary equivalence to $T_{\nu_k}$. Writing $F_kf=(u,v)$, we have \begin{align*}F_kT_{\nu_k}f=&\p{aT_{\nu_{k-1}}u+sT_{\nu_{k-1}}v,se_1 u+de_1 v}\\
=&\tilde{T}_{\nu_k}(u,v)\\
=&\tilde{T}_{\nu_k}F_k f\\
=&\p{T_{\nu_{k-1}}M_{e_{-1}}\oplus I}F_k M_{e_1}f.
\end{align*}
Furthermore, $F_kT_{\nu_k}^N f=\tilde{T}_{\nu_k}^N F_kf.$ 
\end{lem}

\begin{proof}With the help of Lemma \ref{FM_lemma}, we compute:
\begin{align*}&F_kT_{\nu_k}f\\
=&F_k(e_1f-\ip{e_1 f,1}_{L^2(\nu_k)}1)\\
=&F_k(e_1f)-\ip{e_1 f,1}_{L^2(\nu_k)}(1,0)\\
=&\p{ae_1u+se_1v,se_1u+de_1v}-\ip{e_1 f,1}_{L^2(\nu_k)}(1,0)\\
=&\p{ae_1u+se_1v-\ip{e_1 f,1}_{L^2(\nu_k)}1,se_1u+de_1 v}\\
=&\p{ae_1u+se_1v-\ip{F_k e_1 f,F_k 1}_{L^2(\nu_{k-1})\oplus L^2(\nu_{k-1})}1,se_1u+de_1 v}\\
=&\p{ae_1u+se_1v-\ip{ae_1u+se_1v,1}_{L^2(\nu_{k-1})}1,se_1u+de_1 v}\\
=&\p{aT_{\nu_{k-1}}u+sT_{\nu_{k-1}}v,se_1 u+de_1 v}\\
=&\p{T_{\nu_{k-1}}(au+sv),se_1u+de_1 v}\\
=&\p{T_{\nu_{k-1}}(au+sv),M_{e_1}(su+dv)}\\
=&\p{T_{\nu_{k-1}}M_{e_{-1}}\oplus I}F_k M_{e_1}f.
\end{align*}
That $F_kT_{\nu_k}^N=\tilde{T}_{\nu_k}^N F_k$ then follows easily.
\end{proof}

\begin{lem}[human-provided clarification]\label{A_and_S_Lemma}Let $a$ and $s$ be defined as in Lemma \ref{FM_lemma}. Then $\abs{a}\leq 1$ and $\abs{s}=2\sqrt{p_k(1-p_k)}\sin(\pi/Q_k)$.
\end{lem}

\begin{proof}$a=1-p_k(1-e^{2\pi i(1/Q_k)})=1-p_k+p_ke^{2\pi i(1/Q_k)}.$ Therefore, 
\begin{align*}\abs{a}\leq&\abs{1-p_k}+\abs{p_ke^{2\pi i(1/Q_k)}}\\
=&1-p_k+p_k\\
=&1.
\end{align*}
As for $s$, we have $s=\sqrt{p_k(1-p_k)}(e^{2\pi i(1/Q_k)}-1).$ Hence, by the half-angle identities and the fact that $\sin(\pi/Q_k)$ is nonnegative for $Q_k\geq 1$, 
\begin{align*}
    \abs{s}=&\sqrt{p_k(1-p_k)}\abs{e^{2\pi i(1/Q_k)}-1}\\
    =&\sqrt{p_k(1-p_k)}\sqrt{(1-\cos(2\pi/Q_k))^2+\sin^2(2\pi/Q_k)}\\
    =&\sqrt{p_k(1-p_k)}\sqrt{2-2\cos(2\pi/Q_k)}\\
    =&\sqrt{p_k(1-p_k)}\sqrt{4\sin^2(\pi/Q_k)}\\
    =&2\sqrt{p_k(1-p_k)}\sin(\pi/Q_k).
\end{align*}
\end{proof}

\begin{lem}[Claude Opus 5.5 AI-drafted and human-revised]\label{AN_cost_of_one_digit}For any $k\geq1$, \begin{align*}A_N(\nu_k)\leq& A_N(\nu_{k-1})+3\sqrt{p_k}.
\end{align*}
\end{lem}

\begin{proof}Let $F_k e_{-(N+1)}=(u_0,v_0)$. We have 
\begin{align*}u_0=&(1-p_k)e^{-2\pi i(N+1)x}+p_ke^{-2\pi i(N+1)/Q_k}e^{-2\pi i(N+1)x}\\
=&\p{1-p_k+p_ke^{-2\pi i(N+1)/Q_k}}e^{-2\pi i(N+1)x}
\end{align*}
and
\begin{align*}v_0=&\sqrt{p_k(1-p_k)}\p{e^{-2\pi i(N+1)/Q_k}e^{-2\pi i(N+1)x}-e^{-2\pi i(N+1)x}}\\
=&\sqrt{p_k(1-p_k)}\p{e^{-2\pi i(N+1)/Q_k}-1}e^{-2\pi i(N+1)x}.
\end{align*}
Moreover,
\begin{align*}\norm{u_0}=&\abs{1-p_k+p_ke^{-2\pi i(N+1)/Q_k}}\\
\leq& 1
\end{align*}
and
\begin{align*}&\norm{v_0}\\
=&\abs{\sqrt{p_k(1-p_k)}\p{e^{-2\pi i(N+1)/Q_k}-1}}\\
=&2\sqrt{p_k(1-p_k)}\abs{\sin(\pi(N+1)/Q_k)}
\end{align*}
by the argument in the proof of Lemma \ref{A_and_S_Lemma}.

Hence, by Definition \ref{adjointsimplifycor}, Lemma \ref{tilde_T}, and the fact that $T_{\nu_k}$, and therefore $\tilde{T}_{\nu_k}$, are contractions, we have
\begin{align*}A_N(\nu_k)=&\norm{T_{\nu_k}^{N+1}e_{-(N+1)}}_{L^2(\nu_k)}\\
=&\norm{F_k T_{\nu_k}^{N+1}e_{-(N+1)}}_{L^2(\nu_{k-1})\oplus L^2(\nu_{k-1})}\\
=&\norm{\tilde{T}_{\nu_k}^{N+1} F_k e_{-(N+1)}}_{L^2(\nu_{k-1})\oplus L^2(\nu_{k-1})}\\
=&\norm{\tilde{T}_{\nu_k}^{N+1}(u_0,v_0)}_{L^2(\nu_{k-1})\oplus L^2(\nu_{k-1})}\\
\leq&\norm{\tilde{T}_{\nu_k}^{N+1}(u_0,0)}_{L^2(\nu_{k-1})\oplus L^2(\nu_{k-1})}+\norm{\tilde{T}_{\nu_k}^{N+1}(0,v_0)}_{L^2(\nu_{k-1})\oplus L^2(\nu_{k-1})}\\
\leq&\norm{\tilde{T}_{\nu_k}^{N+1}(u_0,0)}_{L^2(\nu_{k-1})\oplus L^2(\nu_{k-1})}+\norm{v_0}_{L^2(\nu_{k-1})}.
\end{align*}

Now letting $D(u,v):=(aT_{\nu_{k-1}}u,de_1 v)$ and $E(u,v):=(sT_{\nu_{k-1}}v,se_1u)$, we have 

$$\tilde{T}_{\nu_k}=D+E.$$

Now clearly

\begin{equation}\label{eq_B}\begin{aligned}&\sum_{j=0}^{n-1}\tilde{T}_{\nu_k}^{n-1-j}ED^j\\
=&\sum_{j=0}^{n-1}\tilde{T}_{\nu_k}^{n-1-j}(\tilde{T}_{\nu_k}-D)D^j\\
=&\sum_{j=0}^{n-1}\p{\tilde{T}_{\nu_k}^{n-j}D^{j}-\tilde{T}_{\nu_k}^{n-1-j}D^{j+1}}\\
=&\tilde{T}_{\nu_k}^{n}-D^{n}
\end{aligned}\end{equation}
because the sum is telescoping. Applying this formula to $(u_0,0)$, we have for any $n$ that
\begin{align*}\tilde{T}_{\nu_k}^{n}(u_0,0)=&D^n(u_0,0)+\sum_{j=0}^{n-1}\tilde{T}_{\nu_k}^{n-1-j}ED^j(u_0,0)\\
=&(a^nT_{\nu_{k-1}}^nu_0,0)+\sum_{j=0}^{n-1}\tilde{T}_{\nu_k}^{n-1-j}E(a^jT_{\nu_{k-1}}^j u_0,0)\\
=&(a^nT_{\nu_{k-1}}^nu_0,0)+\sum_{j=0}^{n-1}\tilde{T}_{\nu_k}^{n-1-j}(0,sa^je_1  T_{\nu_{k-1}}^j u_0).
\end{align*}
Hence, by the Triangle Inequality and the fact that $\norm{u_0}\leq 1$ and $\abs{a}\leq 1$,
\begin{equation}\label{eq_A}\begin{aligned}
    &\norm{\tilde{T}_{\nu_k}^{n}(u_0,0)}\\\leq&\norm{(a^nT_{\nu_{k-1}}^nu_0,0)}+\sum_{j=0}^{n-1}\norm{\tilde{T}_{\nu_k}^{n-1-j}(0,sa^je_1  T_{\nu_{k-1}}^j u_0)}\\
    \leq&\norm{T_{\nu_{k-1}}^n u_0}+\sum_{j=0}^{n-1}\abs{s}\norm{T_{\nu_{k-1}}^{j}}\\
    \leq&\norm{T_{\nu_{k-1}}^n u_0}+2\sqrt{p_k(1-p_k)}\sin(\pi/Q_k)\tau(\nu_{k-1}).
\end{aligned}\end{equation}
Now taking $n=N+1$, and recalling that $u_0$ is $e_{-(N+1)}$ scaled by a number of modulus at most 1, we have that
\begin{align*}\norm{\tilde{T}_{\nu_k}^{N+1}(u_0,0)}\leq&\norm{T_{\nu_{k-1}}^{N+1}e_{-(N+1)}}+2\sqrt{p_k(1-p_k)}\sin(\pi/Q_k)\tau(\nu_{k-1})\\
=&A_N(\nu_{k-1})+2\sqrt{p_k(1-p_k)}\sin(\pi/Q_k)\tau(\nu_{k-1}).
\end{align*}
Therefore, by property $(ii)$ of the sequence $\set{Q_k}$, and using $\abs{\sin(x)}\leq \abs{x}$,

\begin{align*}A_N(\nu_k)\leq& A_N(\nu_{k-1})+2\sqrt{p_k(1-p_k)}\sin(\pi/Q_k)\tau(\nu_{k-1})\\
&+2\sqrt{p_k(1-p_k)}\abs{\sin(\pi(N+1)/Q_k)}\\
\leq&A_N(\nu_{k-1})+2\sqrt{p_k(1-p_k)}\p{\frac{\pi}{Q_k}\cdot\frac{Q_k}{2\pi}+1}\\
\leq& A_N(\nu_{k-1})+3\sqrt{p_k(1-p_k)}\\
\leq& A_N(\nu_{k-1})+3\sqrt{p_k}.
\end{align*}
\end{proof}

We have placed a bound on how much larger $A_N(\nu_k)$ can be than $A_N(\nu_{k-1})$. We now turn to the problem of comparing $A_N(\mu)$ and $A_N(\nu_k)$.

\begin{lem}[Claude Opus 5.5 AI-drafted and human-revised]\label{AN_prod_of_projections}Let $P_n$ denote the orthogonal projection onto $\mathbb{C}e_n$ in $L^2(\mu)$. Then $$A_N(\mu)=\norm{(I-P_0)(I-P_1)\cdots(I-P_N)1}_{L^2(\mu)}.$$
\end{lem}

\begin{proof}Observe that $M_{e_1}^\ast=M_{e_1}^{-1}=M_{e_{-1}}$. Hence,
\begin{align*}M_{e_1}P_kM_{e_1}^{-1}f=&M_{e_1}\ip{M_{e_1}^{-1}f,e_k}e_k\\
=&\ip{f,M_{e_1}e_k}M_{e_1}e_k\\
=&\ip{f,e_{k+1}}e_{k+1}\\
=&P_{k+1}f.
\end{align*}
We have
$T_{\mu}M_{e_1}^{-1}=(I-P_0)M_{e_1}M_{e_1}^{-1}=(I-P_0)$.
By way of induction, suppose $$T_\mu^k M_{e_1}^{-k}=(I-P_0)(I-P_1)\cdots(I-P_{k-1}).$$
Then
\begin{align*}&T_{\mu}^{k+1}M_{e_1}^{-(k+1)}\\
=&(I-P_0)M_{e_1}(I-P_0)(I-P_1)\cdots(I-P_{k-1})M_{e_1}^{-1}\\
=&(I-P_0)M_{e_1}(I-P_0)M_{e_1}^{-1}M_{e_1}(I-P_1)M_{e_1}^{-1}\cdots M_{e_1}(I-P_{k-1})M_{e_1}^{-1}\\
=&(I-P_0)(1-P_1)\cdots(I-P_{k}).
\end{align*}
Consequently,
\begin{align*}A_N(\mu)=&\norm{T_{\mu}^{N+1}e_{-(N+1)}}\\
=&\norm{T_{\mu}^{N+1}M_{e_1}^{-(N+1)}M_{e_1}^{N+1}e_{-(N+1)}}\\
=&\norm{(I-P_0)\cdots(I-P_{N})1}.
\end{align*}
\end{proof}

\begin{lem}[Claude Opus 5.5 AI-drafted and human-revised]\label{AN_Connection}Let $\rho$ be a Borel probability measure on $[0,1)^2$ with singular marginal measures $\rho_x$ and $\rho_y$, that is, $\rho_x(B)=\rho(B\times[0,1))$ and $\rho_y(B)=\rho([0,1)\times B)$. Then for every $N\geq 0$, $$\abs{A_N(\rho_x)-A_N(\rho_y)}\leq 2\sum_{n=0}^{N}\p{\iint\abs{e_n(x)-e_n(y)}^2\,d\rho(x,y)}^{1/2}.$$
\end{lem}

\begin{proof}Clearly, $J_x:L^2(\rho_x)\to L^2(\rho)$ defined by $[J_xf](x,y)=f(x)$ is an isometry, and likewise $J_y:L^2(\rho_y)\to L^2(\rho)$ defined by $[J_y f](x,y)=f(y)$. Define $\mathcal{Q}_{x,n}$ to be the orthogonal projection of $L^2(\rho)$ onto $\mathbb{C}J_x e_n$, and define $\mathcal{Q}_{y,n}$ to be the orthogonal projection of $L^2(\rho)$ onto $\mathbb{C}J_y e_n$. 

For $f\in L^2(\rho_x)$, we have
\begin{align*}\mathcal{Q}_{x,n}J_x f=&\ip{J_xf,J_x e_n}_{L^2(\rho)}J_x e_n\\
=&\ip{f,e_n}_{L^2(\rho_x)}J_x e_n\\
=&J_xP_n f,
\end{align*}
where $P_n$ denotes the orthogonal projection in $L^2(\rho_x)$ onto $\mathbb{C}e_n$. We therefore have $J_x-J_xP_n=J_x-\mathcal{Q}_{x,n}J_x$, which is to say $$J_x(I-P_n)=(I-\mathcal{Q}_{x,n})J_x.$$

Now invoking the result from Lemma \ref{AN_prod_of_projections} and applying the above formula repeatedly, we have

\begin{align*}A_{N}(\rho_x)=&\norm{T_{\rho_x}^{N+1}e_{-(N+1)}}_{L^2(\rho_x)}\\
=&\norm{(I-P_0)(I-P_1)\cdots(I-P_N)1}_{L^2(\rho_x)}\\
=&\norm{J_x(I-P_0)(I-P_1)\cdots(I-P_N)1}_{L^2(\rho)}\\
=&\norm{(I-\mathcal{Q}_{x,0})(I-\mathcal{Q}_{x,1})\cdots(I-\mathcal{Q}_{x,N})J_x1}_{L^2(\rho)}\\
=&\norm{(I-\mathcal{Q}_{x,0})(I-\mathcal{Q}_{x,1})\cdots(I-\mathcal{Q}_{x,N})1}_{L^2(\rho)}.
\end{align*}
By a similar argument,
\begin{align*}A_N(\rho_y)=\norm{(I-\mathcal{Q}_{y,0})(1-\mathcal{Q}_{y,1})\cdots(1-\mathcal{Q}_{y,n})1}_{L^2(\rho)}.
\end{align*}
Now we subtract the operators in these two formulae by regarding them as the simplification of a telescoping sum:
\begin{align*}&\prod_{n=0}^{N}(I-\mathcal{Q}_{x,n})-\prod_{n=0}^{N}(I-\mathcal{Q}_{y,n})\\
=&\sum_{n=0}^{N}\p{\prod_{j=0}^{n}(I-\mathcal{Q}_{x,j})\prod_{j=n+1}^{N}(I-\mathcal{Q}_{y,j})-\prod_{j=0}^{n-1}(I-\mathcal{Q}_{x,j})\prod_{j=n}^{N}(I-\mathcal{Q}_{y,j})}\\
=&\sum_{n=0}^{N}\p{\left[\prod_{j=0}^{n-1}(I-\mathcal{Q}_{x,j})\right]((I-\mathcal{Q}_{x,n})-(I-\mathcal{Q}_{y,n}))\left[\prod_{j=n+1}^{N}(I-\mathcal{Q}_{y,j})\right]}\\
=&\sum_{n=0}^{N}\p{\left[\prod_{j=0}^{n-1}(I-\mathcal{Q}_{x,j})\right](\mathcal{Q}_{y,n}-\mathcal{Q}_{x,n})\left[\prod_{j=n+1}^{N}(I-\mathcal{Q}_{y,j})\right]}.
\end{align*}
Note that the products are operator products that go left-to-right, with the right-most one applied first during evaluation. Because every projection is a contraction, by the Reverse Triangle Inequality, this shows that
\begin{align*}&\abs{A_N(\rho_x)-A_N(\rho_y)}\\
\leq&\norm{\p{\prod_{n=0}^{N}(I-\mathcal{Q}_{x,n})-\prod_{n=0}^{N}(I-\mathcal{Q}_{y,n})}1}\\
\leq&\sum_{n=0}^{N}\norm{\mathcal{Q}_{x,n}-\mathcal{Q}_{y,n}}\\
=&\sum_{n=0}^{N}\norm{\ip{\cdot,J_xe_n}J_xe_n-\ip{\cdot,J_ye_n}J_y e_n}\\
=&\sum_{n=0}^{N}\norm{\ip{\cdot,J_xe_n}(J_xe_n-J_ye_n)+\ip{\cdot,J_x e_n-J_ye_n}J_y e_n}\\
\leq&2\sum_{n=0}^{N}\norm{J_x e_n-J_y e_n}\\
=&\sum_{n=0}^{N}2\p{\iint\abs{e_n(x)-e_n(y)}^2\,d\rho(x,y)}^{1/2}.
\end{align*}
\end{proof}

\begin{lem}[Claude Opus 5.5 AI-drafted and human-revised]\label{AN_difference} Let $k$ be so large that $Q_k\geq N$. Then $$\abs{A_N(\mu)-A_N(\nu_k)}\leq\frac{8\pi}{3(k+1)}.$$
\end{lem}

\begin{proof}
Recall the maps we defined in Definition \ref{with_adjoint_property_measure}. Define a measure $\gamma_k$ on $[0,1)^2$ in the following way: For each Borel $C\subseteq[0,1)^2$, $$\gamma_k(C)=m\p{\set{\omega\in\Omega:\p{\Phi(\omega),\Psi_k(\omega)}\in C}}.$$

Observe that for a Borel $B\subseteq[0,1)$, we have
\begin{align*}\gamma_k(B\times[0,1))=&m\p{\set{\omega\in\Omega:\p{\Phi(\omega),\Psi_k(\omega)}\in B\times[0,1)}}\\
=&m\p{\set{\omega\in\Omega:\Phi(\omega)\in B}}\\
=&m\p{\Phi^{-1}(B)}\\
=&\mu(B).
\end{align*}
Thus $\mu$ is the $x$-marginal measure of $\gamma_k$, and by a similar argument, $\nu_{k}$ is its $y$-marginal measure.

Observe that $m$ is carried by the collection of sequences $X=\set{0,1}^{\infty}$. For any $\omega\in X$, we have

\begin{align*}0\leq&\Phi(\omega)-\Psi_k(\omega)\\
=&\sum_{j=k+1}^{\infty}\frac{\omega_j}{Q_j}\\
\leq&\sum_{j=k+1}^{\infty}\frac{1}{4^{j-k-1}Q_{k+1}}\\
=&\frac{4}{3Q_{k+1}}.
\end{align*}

It follows that $\gamma_k$ is carried by the strip $$S=\set{(x,y):0\leq x-y\leq\frac{4}{3Q_{k+1}}}.$$

On $S$, we have 
\begin{align*}\abs{e_n(x)-e_n(y)}=&\abs{e^{2\pi in(x-y)}-1}\\
\leq&2\pi n\abs{x-y}\\
\leq&\frac{8\pi n}{3Q_{k+1}}.
\end{align*}
Hence, by Lemma \ref{AN_Connection} and property (iv) of Definition \ref{with_adjoint_property_measure}, we have for $Q_k\geq N$,
\begin{align*}\abs{A_N(\mu)-A_N(\nu_k)}\leq&2\sum_{n=1}^{N}\p{\iint\abs{e_n(x)-e_n(y)}^2\,d\gamma_k(x,y)}^{1/2}\\
\leq&2\sum_{n=1}^{N}\frac{8\pi n}{3Q_{k+1}}\\
=&\frac{8\pi N(N+1)}{3Q_{k+1}}\\
\leq&\frac{8\pi Q_{k}(Q_k+1)}{3Q_{k+1}}\\
\leq&\frac{8\pi}{3(k+1)}.
\end{align*}
\end{proof}

We are now prepared to connect our estimates together.

\begin{thm}[Claude Opus 5.5 AI-drafted and human-revised]The measure $\mu$ from Definition \ref{with_adjoint_property_measure} is a singular continuous probability measure on $[0,1)$ that has the adjoint Fourier frame property.
\end{thm}

\begin{proof}Suppose $N$ and $k\geq2$ are positive integers such that $$Q_{k-1}\leq N<Q_{k}.$$ First, observe that by property (iii) of Definition \ref{with_adjoint_property_measure}, and because $Q_{k-1}\leq N$,
\begin{align}\label{AN_elimination}A_N(\nu_{k-2})\leq\norm{T_{\nu_{k-2}}^{N+1}}\leq\norm{T_{\nu_{k-2}}^{Q_{k-1}}}\leq\sqrt{p_{k-1}}.\end{align}

Because $N\leq Q_k$, Lemma \ref{AN_difference} says $$\abs{A_N(\mu)-A_N(\nu_k)}\leq\frac{8\pi}{3(k+1)}.$$ 

Finally, Lemma \ref{AN_cost_of_one_digit} says for any $j\geq1$, \begin{align*}A_N(\nu_j)\leq& A_N(\nu_{j-1})+3\sqrt{p_j}.
\end{align*} Combining these three results, we have that \begin{align*}A_N(\mu)\leq&A_N(\nu_k)+\frac{8\pi}{3(k+1)}\\
\leq&A_N(\nu_{k-1})+3\sqrt{p_k}+\frac{8\pi}{3(k+1)}\\
\leq&A_N(\nu_{k-2})+3\sqrt{p_{k-1}}+3\sqrt{p_k}+\frac{8\pi}{3(k+1)}\\
\leq&4\sqrt{p_{k-1}}+3\sqrt{p_k}+\frac{8\pi}{3(k+1)}.
\end{align*}
For every integer $N\geq Q_{1}$ there is a unique $k_N$ such that $Q_{k_{N}-1}\leq N<Q_{k_{N}}$, and clearly $\lim_{N\to\infty}k_N=\infty$. Therefore,
$\lim_{N\to\infty}A_N(\mu)=0$ by the fact that $\lim_{k\to\infty} p_{k}=0$.

Hence, $\mu$ has the adjoint Fourier frame property.
\end{proof}

\subsection{A Singular Continuous Measure Without the Adjoint Fourier Frame Property}\index{Adjoint frame property for stationary effective sequences}

We now show by way of example that not all singular continuous measures have the adjoint Fourier frame property. We do this by modifying the previous construction, the principal difference being that the probability weights $p_k$ no longer approach $0$.

\begin{defn}[Claude Opus 5.5 AI-drafted and human-revised]\label{without_adjoint_prop}
Construct a new measure $\mu$ and a sequence of finite atomic measures $\set{\nu_k}$ exactly as in Definition \ref{with_adjoint_property_measure}, except with the following differences:

\begin{enumerate}
   \item $p_k=\frac{1}{2}$ for all $k$ (thus $p_k\not\to 0)$.
   \item $Q_1$ is even, and hence all $Q_k$ are even for $k\geq 1$.
   \item Property (ii) is strengthened to $Q_k\geq 2\pi k\tau(\nu_{k-1})$ for $k\geq 1$.
   \item Property (iii) may be dropped.
\end{enumerate}
\end{defn}

\begin{lem}The measure $\mu$ of Definition \ref{without_adjoint_prop} is a singular continuous probability measure.
\end{lem}
\begin{proof}Lemma \ref{singular_continuous} applies without change.
\end{proof}

\begin{lem}[Claude Opus 5.5 AI-drafted and human-revised]\label{lower_bound_on_doubling_atomic}Define $H_k:L^2(\nu_{k-1})\to L^2(\nu_{k})$ as follows: For $f\in L^2(\nu_{k-1})$, set
\begin{align*}[H_kf](x):=-f(x)\\
[H_kf](x+1/Q_k):=f(x)
\end{align*}
for $x\in\text{supp}(\nu_{k-1}).$ Then for every $n\geq 0$, $$\norm{T_{\nu_{k}}^nH_k f}_{L^2(\nu_{k})}\geq\p{\abs{\cos(\pi/Q_k)}^n-\sin(\pi/Q_k)\tau(\nu_{k-1})(1+n\sin(\pi/Q_k))}\norm{f}_{L^2(\nu_{k-1})}.$$
\end{lem}

\begin{proof}Assume the setup of Lemma \ref{Fk_Operator}. Since $p_k=\frac{1}{2}$, we have 
\begin{equation}\label{FkHk_eq}\begin{aligned}&F_kH_kf\\
=&\p{\frac{1}{2}(H_k f(x)+H_k f(x+1/Q_k)),\frac{1}{2}(H_k f(x+1/Q_k)-H_k f(x))}\\
=&\p{0,f}.\end{aligned}\end{equation}
Assume also the setup of Lemma \ref{FM_lemma}, where $\omega=e^{2\pi i/Q_k}$. We have
\begin{align*}a=&\frac{1}{2}(1+\omega)\\
s=&\frac{1}{2}(\omega-1)\\
d=&\frac{1}{2}(1+\omega).
\end{align*}
Then $\abs{a}=\abs{d}=\abs{\cos(\pi/Q_k)}$ and $\abs{s}=\abs{\sin(\pi/Q_k)}=\sin(\pi/Q_k)$. As in Lemma \ref{AN_cost_of_one_digit}, we may write $\tilde{T}_{\nu_k}=D+E$, where $D(u,v)=(aT_{\nu_{k-1}}u,a e_1 v)$ and $E(u,v)=(sT_{\nu_{k-1}}v,se_1 u)$.

Now observe that using equation \eqref{eq_A} with $u_0$ replaced by $w$, we have

\begin{equation}\label{eq_C}\norm{\tilde{T}_{\nu_k}^m(w,0)}\leq\norm{T_{\nu_{k-1}}^m w}+\abs{s}\tau(\nu_{k-1})\norm{w}.\end{equation}

Now by Equation \eqref{eq_B}, but this time applied to a vector $(0,v)$, we have $$\tilde{T}_{\nu_{k}}^n(0,v)=D^n(0,v)+\sum_{j=0}^{n-1}\tilde{T}_{\nu_{k}}^{n-1-j}ED^{j}(0,v).$$
Since $D^j(0,v)=(0,a^je_jv)$ and $E(0,v)=(s T_{\nu_{k-1}}v,0)$, it follows that $ED^j(0,v)=(sa^j T_{\nu_{k-1}}e_j v,0)$.

Applying Inequality \eqref{eq_C} to this,
\begin{align*}\norm{\tilde{T}_{\nu_k}^{n}(0,v)-D^n(0,v)}\leq&\sum_{j=0}^{n-1}\norm{\tilde{T}_{\nu_{k}}^{n-1-j}ED^{j}(0,v)}\\
=&\sum_{j=0}^{n-1}\norm{\tilde{T}_{\nu_k}^{n-1-j}(sa^j T_{\nu_{k-1}}e_jv,0)}\\
\leq&\sum_{j=0}^{n-1}\p{\norm{T_{\nu_{k-1}}^{n-1-j}sa^jT_{\nu_{k-1}}e_jv}+\abs{s}\tau(\nu_{k-1})\norm{sa^jT_{\nu_{k-1}}e_jv}}\\
\leq&\sum_{j=0}^{n-1}\p{\abs{s}\norm{T_{\nu_{k-1}}^{n-1-j}}\norm{v}+\abs{s}^2\tau(\nu_{k-1})\norm{v}}\\
\leq&\abs{s}\norm{v}\tau(\nu_{k-1})+n\abs{s}^2\tau(\nu_{k-1})\norm{v}\\
=&\abs{s}\tau(\nu_{k-1})\norm{v}\p{1+n\abs{s}}.
\end{align*}
Finally, this shows that 
\begin{equation}\label{eq_D}\begin{aligned}\norm{\tilde{T}_{\nu_k}^n(0,v)}\geq&\norm{D^n(0,v)}-\abs{s}\tau(\nu_{k-1})\norm{v}\p{1+n\abs{s}}\\
=&\abs{a}^n\norm{v}-\abs{s}\tau(\nu_{k-1})\norm{v}\p{1+n\abs{s}}.
\end{aligned}\end{equation}
We note that by Lemma \ref{tilde_T} and equation \eqref{FkHk_eq},
$$F_k T_{\nu_k}^{n}H_kf=\tilde{T}_{\nu_k}^n F_k H_k f=\tilde{T}_{\nu_k}^n (0,f).$$
Therefore, because $F_k$ is unitary, we have $$\norm{\tilde{T}_{\nu_k}^n(0,f)}_{L^2(\nu_{k-1})\oplus L^2(\nu_{k-1})}=\norm{T_{\nu_k}^nH_kf}_{L^2(\nu_k)}.$$
Thus, inequality \eqref{eq_D} is equivalent to the assertion of the lemma after substituting $v=f$, $\abs{a}=\abs{\cos(\pi/Q_k)}$, and $\abs{s}=\sin(\pi/Q_k)$.
\end{proof}

\begin{thm}[Claude Opus 5.5 AI-drafted and human-revised]The measure $\mu$ of Definition \ref{without_adjoint_prop} is a singular continuous probability measure on $[0,1)$ that does not have the adjoint Fourier frame property.
\end{thm}

\begin{proof}Since each $Q_k$ is even, we may define integers $N_k$ such that $\frac{N_k+1}{Q_k}=\frac{1}{2}$. It is then not difficult to see that

$$H_k(-e_{-(N_k+1)})(x)=e^{-2\pi i(N_k+1)x}$$ and 
\begin{align*}H_k(-e_{-(N_k+1)})(x+1/Q_k)=&-e^{-2\pi i(N_k+1)x}\\
=&e^{-2\pi i(1/2)}e^{-2\pi i(N_k+1)x}\\
=&e^{-2\pi i(N_k+1)/Q_k}e^{-2\pi i(N_k+1)x}\\
=&e^{-2\pi i(N_k+1)(x+1/Q_k)}.
\end{align*}
In other words, $H_k(-e_{-(N_k+1)})=e_{-(N_k+1)}.$ Then by Lemma \ref{lower_bound_on_doubling_atomic} and the new property (ii), we obtain
\begin{align*}A_{N_k}(\nu_k)=&\norm{T_{\nu_k}^{N_k+1}e_{-(N_k+1)}}\\
=&\norm{T_{\nu_k}^{N_k+1}H_k(-e_{-(N_k+1)})}\\
\geq&\abs{\cos(\pi/Q_k)}^{N_k+1}-\sin(\pi/Q_k)\tau(\nu_{k-1})(1+(N_k+1)\sin(\pi/Q_k))\\
=&\abs{\cos(\pi/Q_k)}^{Q_k/2}-\sin(\pi/Q_k)\tau(\nu_{k-1})\p{1+\frac{Q_k}{2}\sin(\pi/Q_k)}\\
\geq&(\cos(\pi/Q_k))^{Q_k/2}-\frac{\pi}{Q_k}\tau(\nu_{k-1})\p{1+\frac{Q_k}{2}\cdot\frac{\pi}{Q_k}}\\
=&(\cos(\pi/Q_k))^{Q_k/2}-\tau(\nu_{k-1})\p{\frac{\pi}{Q_k}+\frac{\pi^2}{2Q_k}}\\
\geq&(\cos(\pi/Q_k))^{Q_k/2}-\frac{Q_k}{2\pi k}\p{\frac{\pi}{Q_k}+\frac{\pi^2}{2Q_k}}\\
=&(\cos(\pi/Q_k))^{Q_k/2}-\frac{1}{2k}-\frac{\pi}{4k}.
\end{align*}
 A version of Bernoulli's Inequality states that $(1+x)^r\geq 1+rx$ for every integer $r\geq1$ and every real number $x\geq-1$. Since $\cos(x)\geq1-\frac{x^2}{2}$, we have that
 \begin{align*}(\cos(\pi/Q_k))^{Q_k/2}\geq&\p{1-\frac{\pi^2}{2Q_k^2}}^{Q_k/2}\\
 \geq&1-\frac{Q_k}{2}\frac{\pi^2}{2Q_k^2}\\
 =&1-\frac{\pi^2}{4Q_k}.
 \end{align*}
 Thus,

 $$A_{N_k}(\nu_k)\geq1-\frac{\pi^2}{4Q_k}-\frac{1}{2k}-\frac{\pi}{4k}.$$

 The argument in the proof of Lemma \ref{AN_difference} still shows that
 \begin{align*}\abs{A_{N_k}(\mu)-A_{N_k}(\nu_k)}\leq&\frac{8\pi N_k(N_k+1)}{3Q_{k+1}}\\
 =&\frac{8\pi\cdot\p{\frac{Q_k}{2}-1}\frac{Q_k}{2}}{3Q_{k+1}}\\
 \leq&\frac{2\pi Q_k^2}{3Q_{k+1}}\\
 \leq&\frac{2\pi Q_k^2}{3(k+1)Q_k(Q_k+1)}\\
 \leq&\frac{2\pi}{3(k+1)}.
 \end{align*}
 Combining these two inequalities,
 \begin{align*}A_{N_k}(\mu)\geq&A_{N_k}(\nu_k)-\frac{2\pi}{3(k+1)}\\
 \geq&1-\frac{\pi^2}{4Q_k}-\frac{1}{2k}-\frac{\pi}{4k}-\frac{2\pi}{3(k+1)}.
 \end{align*}
 Since $A_N(\mu)=\norm{T_{\mu}^{N+1}e_{-(N+1)}}\leq 1$ for all $N$, $\lim_{k\to\infty}A_{N_k}(\mu)=1$. In particular, $\lim_{N\to\infty}A_N(\mu)\neq 0$. Thus, $\mu$ does not have the adjoint Fourier frame property.
\end{proof}

\begin{rem}We have $\limsup_{N\to\infty}A_N(\mu)=1$, which is the ``worst'' this limit superior can be.
\end{rem}